\documentclass[12pt]{amsart}
\usepackage{amssymb,amsmath,amscd}
\newtheorem{thm}{Theorem}[section]
\newtheorem{corollary}[thm]{Corollary}
\newtheorem{prop}[thm]{Proposition}
\newtheorem{lemma}[thm]{Lemma}
\newtheorem{fact}[thm]{Fact}
\newtheorem{conj}[thm]{Conjecture}
\newtheorem{question}[thm]{Question}

\newtheorem{defn}[thm]{Definition}
\newtheorem{example}[thm]{Example}

\newtheorem{remark}[thm]{Remark}

\newcommand{\bt}{\begin{thm}}
\newcommand{\et}{\end{thm}}
\newcommand{\bp}{\begin{prop}}
\newcommand{\ep}{\end{prop}}
\newcommand{\bd}{\begin{defn}}
\newcommand{\ed}{\end{defn}}
\newcommand{\bl}{\begin{lemma}}
\newcommand{\el}{\end{lemma}}
\newcommand{\bfa}{\begin{fact}}
\newcommand{\efa}{\end{fact}}
\newcommand{\bc}{\begin{corollary}}
\newcommand{\ec}{\end{corollary}}
\newcommand{\bex}{\begin{example}}
\newcommand{\eex}{\end{example}}
\newcommand{\br}{\begin{remark}}
\newcommand{\er}{\end{remark}}
\newcommand{\bq}{\begin{question}}
\newcommand{\eq}{\end{question}}
\newcommand{\bconj}{\begin{conj}}
\newcommand{\econj}{\end{conj}}

\newcommand{\Hom}{{\mathcal H}om}

\newcommand{\supp}{{\mbox{Supp}}}
\newcommand{\rank}{{\mbox{rank}}}
\newcommand{\Pic}{{\mbox{Pic}}}

\newcommand{\Length}{{\mbox{length}}}

\newcommand{\Ext}{{\mathcal E}xt}

\newcommand{\sotto}[2]{#1_{#2}}
\newcommand{\lra}{\longrightarrow}
\newcommand{\rrr}{\rightarrow}

\newcommand{\ideal}[1]{\sotto {{\mathcal I}}{#1}}
\newcommand{\idealc}{\ideal{C}}

\newcommand{\exact}[3]
{0 \rrr #1 \rrr #2 \rrr #3 \rrr 0}

\newcommand{\pso}{\mathbb{P}^3}
\newcommand{\pthree}{\mathbb{P}^3}
\newcommand{\ptwo}{\mathbb{P}^2}
\newcommand{\pone}{\mathbb{P}^1}
\newcommand{\PP}{\mathbb{P}}

\newcommand{\ra}{\rightarrow}

\newcommand{\psn}{\mathbb{P}^N}
\newcommand{\Pone}{\mathbb{P}^1}
\newcommand{\Aone}{\mathbb{A}^1}

\newcommand{\Pthree}{\mathbb{P}^3}

\newcommand{\Z}{\mathbb{Z}}

\newcommand{\ds}{\displaystyle}

\newcommand{\B}{{\mathcal B}}

\newcommand{\coo}{{\mathcal O}}

\newcommand{\caf}{{\mathcal F}}

\newcommand{\mathcall}{{\mathcal L}}
\newcommand{\mathcalz}{{\mathcal Z}}
\newcommand{\mathcaly}{{\mathcal Y}}
\newcommand{\mathcals}{{\mathcal S}}
\newcommand{\mathcalw}{{\mathcal W}}

\newcommand{\cae}{{\mathcal E}}
\newcommand{\cag}{{\mathcal G}}

\begin{document}

\title{Hilbert Schemes of Degree Four Curves }
\author{Scott Nollet}
\address{Department of Mathematics,
Texas Christian University, Fort Worth, TX 76129, USA}
\email{s.nollet@tcu.edu}

\author{Enrico Schlesinger}
\address{Dipartimento di Matematica, Politecnico di Milano, Piazza
Leonardo da Vinci 32,
20133 Milano, Italy}
\email{enrsch@mate.polimi.it}

\date{}

\begin{abstract}
          In this paper we determine the irreducible components of the
          Hilbert schemes $H_{4,g}$ of locally Cohen-Macaulay space curves
          of degree four and arbitrary arithmetic genus $g$.
          We show that these Hilbert schemes are connected, in spite
          of having $ \sim \frac{g^{2}}{24}$ irreducible components.
          For $g \leq -3$ we exhibit a component that is disjoint from
          the component of extremal curves and use this to give
          a counterexample to a conjecture of  A\"it-Amrane and Perrin.
\end{abstract}

\maketitle

\section{Introduction} \label{intro}

Liaison theory has played a prominent role in the
classification of algebraic curves in $\pso_k$ since the pioneering
work of Max Noether. It has only recently become clear that
locally Cohen-Macaulay curves - locally Cohen-Macaulay schemes of
pure dimension one - are the natural object of study ~\cite{MDP},
even if one is only interested in smooth connected curves.
For example, in Gruson and Peskine's classification of
smooth irreducible curves of degree 8 and genus 5 \cite{gp}, one family
of curves is in the biliaison class of double lines of genus $-2$
(non-reduced curves) and another is in the biliaison class of the disjoint
unions of a line and a twisted cubic (non-connected curves).
In general, every biliaison class contains
an essentially unique  minimal curve, from which every  other curve in
the class
maybe obtained by a rather explicit procedure - known to the expert
as the Lazarsfeld-Rao property \cite{bbm}.
While every biliaison class contains smooth
connected curves, minimal curves need only be locally Cohen-Macaulay.
This explains the interest in the Hilbert schemes $H_{d,g}$ parametrizing
locally Cohen-Macaulay curves in $\pso$ of degree $d$ and
{\em arithmetic genus} $g$.

Here are some results on these Hilbert schemes.
In general, $H_{d,g}$ is non-empty precisely when $g=\frac{1}{2}(d-1)(d-2)$
(plane curves) or $d > 1$ and $g \leq \frac{1}{2} (d-2)(d-3)$
\cite[3.3 and 3.4]{genus}. It is also known that $H_{d,g}$ is
reducible for $d \geq 3$ and $g \leq \frac{1}{2}(d-3)(d-4)$ with
two exceptions: $(d,g)=(3,0)$ and $(3,-1)$ \cite{extremal}.
More recently it has been shown that $H_{3,g}$ is connected
(it has $\lfloor \frac{4-g}{3} \rfloor$ irreducible components) \cite{nthree}
and that $H_{d,g}$ is connected for
$g > \binom{d-3}{2}-2$ \cite{nolletp,phsamir,sabadini}.

Each connectedness result above was obtained by specializing various
families of curves to {\it extremal curves} as introduced by
Martin-Deschamps and Perrin \cite{bounds}: these are the
curves $C$ which have the largest Rao function $h^1 \ideal{C} (n)$
with respect to $d$ and $g$. The extremal curves are geometrically
characterized  as the  curves  of  degree $d$  that  contain a  planar
subcurve of degree
$d-1$ - unless $g= \binom{d-1}{2}$, when $C$
itself is planar, or $(d,g) \in \{(3,0),(4,1)\}$-
and their closure forms an irreducible component $E \subset H_{d,g}$
\cite{extremal}. The existence of a component of curves with the largest
Rao function led Hartshorne to ask the following questions: \\

{\bf Question 1:} Is $H_{d,g}$ connected for all $d$ and $g$ ? \\

{\bf Question 2:} Does each component
$B \subset H_{d,g}$ meet the extremal component $E$ ? \\

Hartshorne showed that various
families of curves can be connected to extremal curves, for example
smooth rational and elliptic curves, smooth curves of degree
$d \geq g+3$, arithmetically Cohen-Macaulay curves and many others
\cite{hconn}.
Related to this is a conjecture of A{\"\i}t-Amrane and Perrin stating
that if $X$ is a family of curves whose cohomology does not exceed
that of a family $X_{0}$ and the Rao module of the general curve in
$X$ is a flat deformation of a subquotient of the Rao module of curves
in $X_{0}$, then ${\overline X} \cap X_{0} \neq \emptyset$ (they have
shown \cite{aitperrin} that semi-continuity alone is insufficient).

Now we specialize to curves of degree $d=4$. The Hilbert scheme
$H_{4,3}$ parametrizes plane curves and is smooth irreducible of
dimension $17$; $H_{4,1}$ is smooth irreducible of dimension $16$ by
\cite{Ell} and \cite[3.3 and 3.5]{genus} and its general curve is a complete
intersection of two quadrics.
The Hilbert scheme $H_{4,0}$ has two irreducible components, whose general
members are respectively rational quartic curves and disjoint unions of a
plane cubic and a line.
Hartshorne first noticed that these two families can be connected;
now there are several published proofs: \cite[3.10]{nthree},
\cite[5.21 and 5.22]{hmdp3}, \cite[1.1]{hconn}.
One of our motivations is to extend the systematic study of
these Hilbert schemes and complete the picture when $d=4$.

Our first theorem is the classification of curves in $H_{4,g}$: we
describe the irreducible components and give their dimensions.
Our method is to first identify components whose general curve
$C$ is very special in the sense that
$C$ is contained in quadric surface ($h^{0} (\ideal{C} (2)) \geq 1$)
or $C$ has large speciality ($h^{1} \coo_{C} (-1) \geq 2$).
There are very few such components: if the
general curve $C$ of a family lies on a quadric surface, then either
$C$ is an {\it extremal} curve, a {\it subextremal} curve (these are
the curves with largest Rao function among the non-extremal curves
\cite{subextremal}) or a double conic.

On the other hand, if the general curve $C$ of an irreducible
component of $H_{4,g}$ has large speciality but does not lie on a
quadric, we show (Proposition \ref{d3g0-c}) that $C$ is either a {\it thick
$4$-line} or the union of a conic and a double line meeting at a
point with multiplicity $2$. A thick $4$-line~\cite{banica} is a curve
of degree $4$ supported on a line $L$ and containing the first infinitesimal
neighborhood of $L$ in $\pso$.

Having disposed of these few very special components, it is relatively
easy to list the other irreducible components. Their general member
is either (a) a quasiprimitive (i.e. non-thick~\cite{banica}) $4$-line,
(b) the disjoint union of a line and a general curve of an irreducible
component of $H_{3,g}$, or (c) the disjoint union of two double lines.
Asymptotically the Hilbert scheme $H_{4,g}$ has $\sim
\frac{g^{2}}{24}$ irreducible components, most of which are families
of 4-lines (there are roughly $\frac{-3g}{2}$ components which
are not). For example, from our table (Theorem \ref{comps}) we find that
$H_{4,-100}$ has $530$ components ($377$ of these arise from 4-lines)
while $H_{4,-1000}$ has $42755$ components
(of which $41252$ arise from 4-lines).

Our second theorem states that $H_{4,g}$ is connected whenever it is
not empty. For $g \leq -3$, the main novelty is the presence of a
component $G_4$ of thick $4$-lines that consists entirely of curves
with generic embedding dimension three. We prove the connectedness
theorem by showing that each irreducible component can be connected
either to the extremal component $E$ or to $G_{4}$, and that the
component of subextremal
curves meets both $E$ and $G_4$. Specifically, the quasi-primitive
4-lines and the
curves with large speciality may be deformed to thick 4-lines
($\S$ \ref{4lines} and $\S$ \ref{twisted}). We show that a disjoint union of
double lines specializes to a quasi-primitive 4-line on a double
quadric surface in section \ref{sdouble} and section \ref{striple}
is devoted to showing that families of unions of triple lines and
reduced lines can be connected to the extremal component.

The component $G_{4}$ of thick 4-lines turns out to be rather
interesting. Since these curves are
scheme-theoretically (although not cohomologically) the most special,
they cannot specialize to extremal curves, answering Question 2 in the
negative.
Since their Rao modules are flat deformations of subquotients of the
Rao modules of extremal curves (Example \ref{perrin}), we obtain a
counterexample to the conjecture of A{\"\i}t-Amrane and Perrin above.
Question 1 remains open. \\

\noindent{\bf Notation and conventions} \newline We work over an
algebraically closed field $k$ of arbitrary characteristic. A
curve for us is a locally Cohen-Macaulay scheme over $k$ of pure
dimension $1$.

We will freely use the sentence ``the family of
curves of degree $d$ and genus $g$ with property $P$ is
irreducible of dimension $m$'' meaning there is a (unique)
irreducible $m$-dimensional constructible subset $S$ of the
Hilbert scheme $H_{d,g}$ whose closed points parametrize the
curves of degree $d$ and genus $g$ with property $P$. Note that,
since $S$ is constructible, the closure of $S$ in the Hilbert
scheme is also an $m$-dimensional irreducible subset.

The symbol $L \cup_{nP} C$ denotes the schematic union of a
line $L$ and a curve $C$,
whose intersection is the divisor $nP$ on $L$.

\section{Multiplicity four structures on lines} \label{4lines}

In this section we study locally Cohen-Macaulay curves in
$\Pthree$ which are supported on a line; we will simply call these
d-lines, where $d$ is the degree of the curve.

We begin with the general theory of Banica and Forster \cite[$\S
\; 3$]{banica}. Let $C$ be a locally Cohen-Macaulay curve on a
smooth threefold $X$ with smooth support $Y$. Letting $Y^{(i)}$ be
the subscheme of $X$ defined by $\ideal{Y}^{i}$ and $C_{i}$ be the
subscheme of $X$ obtained by removing the embedded points from $C
\cap Y^{(i)}$, we obtain the {\it Cohen-Macaulay filtration} for
$Y \subset C$: $$Y = C_{1} \subset C_{2} \subset \dots \subset
C_{k} = C$$ for some $k \geq 1$. The quotients $L_{j} =
\ideal{C_{j}}/\ideal{C_{j+1}}$ are vector bundles on $Y$ and the
multiplicity is given by $\mu(C) = 1 + \sum \rank L_{j}$. The
natural inclusions $\ideal{C_{i}} \ideal{C_{j}} \subset
\ideal{C_{i+j}}$ induce generically surjective maps $L_{i} \otimes
L_{j} \ra L_{i+j}$ and hence we obtain generic surjections
$L_{1}^{j} \ra L_{j}$.

As in \cite[$\S 4$]{banica}, we say that $C$ is {\it thick} if it
contains $Y^{(2)}$, i.e. $C_{2} = Y^{(2)}$. This is also
equivalent to the condition that $C$ have embedding dimension
three at each point. In this case $L_{1} =
\ideal{Y}/{\ideal{Y}^{2}}$ is the conormal bundle of $Y$ on $X$.
If further $\mu(C)=4$, then $\rank L_{2} = 1$ and there is an
exact sequence $$0 \lra \frac{\ideal{C}}{\ideal{Y}^{3}} \lra
\frac{\ideal{Y}^{2}}{\ideal{Y}^{3}} \lra L_{2} \lra 0.$$ If $Y$ is
a line in $\Pthree$, we obtain the following.

\bp\label{thickmoduli} For $g \leq 1$ the set of thick $4$-lines
of genus $g$ is parametrized by an irreducible closed subscheme of
$H_{4,g}$ of dimension $9-3g$.
%For $g \leq -3$ this subscheme is an irreducible component of $H_{4,g}$.
\ep
\begin{proof}
The condition that $C$ contain $Y^{(2)}$ is clearly closed. If $C$
is a thick 4-line with support $Y \subset \Pthree$, then
$\mu(C)=4$, $\frac{\ideal{Y}^{2}}{\ideal{Y}^{3}} \cong
\coo_{Y}(-2)^{\oplus 3}$ and $L_{2} \cong \coo_{Y}(-g-1)$ (because
$\Pic \; Y \cong \Z$). It follows from the exact sequence above
that giving such a curve $C$ is the same as giving a surjective
morphism $\coo_Y(-2)^{\oplus 3} \rrr \coo_Y (-g-1)$ modulo an
automorphism of $\coo_Y (-g-1)$. Thus the set of thick $4$-lines
of genus $g$ is parametrized by an open subset of a
$\mathbb{P}^{5-3g}$-bundle over the Hilbert scheme of lines in
$\pso$.
\end{proof}

If $C$ has generic embedding dimension two it is said to be {\it
quasiprimitive}. In this case $\rank L_{1} = 1$ and the generic
surjections $L_{1}^{j} \ra L_{j}$ yield effective divisors $D_{j}$
such that $L_{j} \cong L_{1}^{j}(D_{j})$; the multiplication maps
show that $D_{i} + D_{j} \leq D_{i+j}$. If $C$ is a quasiprimitive
4-line in $\pso$, then the Cohen-Macaulay filtration takes the
form $Y \subset D \subset W \subset C$ where $Y$ is a line, $D$ a
double line and $W$ a triple line. Setting $a=\deg L_{1}, b = \deg
D_{2}$ and $c = \deg D_{3}$, we call the triple $(a,b,c)$ the {\it
type} of $C$. Note that $a \geq -1$ because the surjection
$\ideal{Y} \ra \coo_{Y}(a)$ factors through
$\ideal{Y}/\ideal{Y}^{2} \cong \coo_{Y}(-1)^{2}$. The isomorphisms
$\ideal{Y,D} \cong \coo_Y (a)$, $\ideal{D,W} \cong \coo_Y (2a+b)$
and $\ideal{W,C} \cong \coo_Y (3a+c)$ show that $p_{a}(C)= -6a -b
-c -3$.

If $C$ has type $(a,b,c)$ with $a=-1$, then $C$ necessarily lies
in a double plane and hence is a flat limit of double conics
\cite[8.1 and 8.2]{2h}. We are mainly interested in families that
form irreducible components of the Hilbert scheme, so we will
assume that $a \geq 0$ in the sequel. The proof of the following
proposition is based on~\cite[\S 3.8]{banica}.

\bp\label{qp} For a triple of integers $(a,b,c)$ satisfying $a
\geq 0$ and $0 \leq b \leq c$, the set of quasiprimitive $4$-lines
of type $(a,b,c)$ is parametrized by an irreducible constructible
subscheme of $H_{4,g}$ of dimension $9a+2b+2c+13$. \ep

\begin{proof}
The set of double lines of type $a$ is parametrized by an open
subscheme $V_2$ of a $\mathbb{P}^{2a+3}$-bundle over the
Grassmannian of lines in $\pso$~\cite[1.6]{nthree}. Indeed, to
give a double structure $D$ of type $a$ on a line $Y$ is
equivalent to give a surjective morphism $\ideal{Y}/\ideal{Y}^2
\ra \coo_Y (a)$ modulo an automorphism of $\coo_{Y} (a)$.

Similarly, the set of quasiprimitive triple lines of type $(a,b)$
containing a double line $D \in V_{2}$ is given by surjective
morphisms $\ideal{D}/ {\ideal{Y} \ideal{D}} \ra \coo_{Y}(2a+b)$
modulo automorphisms of $\coo_{Y}(2a+b)$. Since $\ideal{D}/
{\ideal{Y} \ideal{D}} \cong \coo_Y (2a) \oplus \coo_Y(-a-2)$
\cite[2.3 and 2.6]{nthree}, it follows that the set of triple
lines of type $(a,b)$ is parametrized by an open subscheme $V_3$
of a $\mathbb{P}^{3a+2b+3}$-bundle over $V_2$, hence is
irreducible of dimension $5a + 2b +10$.

In the same way, 4-lines of type $(a,b,c)$ containing a fixed
quasiprimitive 3-line $W$ of type $(a,b)$ with support $Y$ are
given by surjections $\nu_{W}=\ideal{W}/{\ideal{Y} \ideal{W}}
\stackrel{\phi}{\ra} \coo_{Y}(3a+c)$ modulo automorphisms of
$\coo_{Y}(3a+c)$. Noting that $\ideal{D}^{2} \subset {\ideal{Y}
\ideal{W}}$ on an open set, we see that the image of
$\ideal{D}^{2}$ in $\nu_{W}$ is torsion, hence $\phi$ factors
through $\ideal{W}/\mathcal{J}$, where
$\mathcal{J}=\ideal{Y}\ideal{W} + \ideal{D}^{2}$. We claim that
(a) $\ideal{W}/\mathcal{J} \cong \coo_{Y}(3a+b) \oplus
\coo_{Y}(-a-b-2)$ and (b) the induced map $\ideal{Y}^{3} \ra
\coo_{Y}(-a-b-2)$ is zero. Claim (b) shows that each surjection
$\psi$ defines a {\it quasiprimitive} 4-line $C$ (since
$\coo_{Y}(-a-b-2) \ra \coo_{Y}(3a+c)$ cannot be surjective for $a
\geq 0$) and (a) shows that the set of surjections is parametrized
by an open subscheme $V_{4}$ of a $\mathbb{P}^{4a+2c+3}$-bundle
over $V_3$, hence is irreducible of dimension $9a + 2b +2c +13$.

Noting the exact sequence
\begin{equation}\label{splitter}
      {0 \ra \frac{\ideal{Y} \ideal{D}}{\mathcal{J}}
      \ra \frac{\ideal{W}}{\mathcal{J}} \ra
      \frac{\ideal{W}}{\ideal{Y} \ideal{D}}  \ra 0,}
\end{equation}
we observe the following: The multiplication map $\ideal{Y,D}
\otimes_{\coo_{Y}} \ideal{D,W} \ra {{\ideal{Y} \ideal{D}}/
\mathcal{J}}$ is an isomorphism, since it is surjective and
${{\ideal{Y} \ideal{D}}/ \mathcal{J}}$ has rank one (recall that
$\mathcal{J}=\ideal{Y}\ideal{W}$ on an open set) and hence the
kernel is zero. It follows that ${{\ideal{Y} \ideal{D}}/
\mathcal{J}} \cong \coo_{Y}(3a+b)$. Further, the exact sequence $$
\exact {  \ideal{W}/{\ideal{Y} \ideal{D}} } { \ideal{D}/{\ideal{Y}
\ideal{D}}  } {  \coo_{Y} (2a+b) } $$ shows that ${
\ideal{W}/{\ideal{Y} \ideal{D}} } \cong \coo_{Y}(-a-b-2)$, since
$\ideal{D}/{\ideal{Y} \ideal{D}} \cong \coo_{Y} (2a) \oplus
\coo_{Y} (-a-2)$. Finally, the exact sequence \ref{splitter}
splits because $a,b \geq 0$, which proves part (a) of the claim.
Part (b) is clear because $\ideal{Y}^{3} \subset
\ideal{Y}\ideal{D}$.

\end{proof}

We now show that the irreducible family of quasiprimitive
$4$-lines of type $(a,b,c)$ contains thick $4$-lines in its
closure for $a \geq 0$. This proposition will ultimately allow us
to prove the connectedness of the Hilbert scheme $H_{4,g}$.

\bp \label{thintothick} Let $(a, b, c)$ be a triple satisfying $a
\geq 0$ and $c \geq b \geq 0$. Then there exists a flat family of
curves $C \subset \pso_{\Aone}$ such that
\begin{enumerate}
\item  the fibre  $C_{t}$ is a quasiprimitive 4-line of type
$(a,b,c)$ for $t \neq 0$ and
\item the fibre $C_{0}$ is a thick  4-line.
\end{enumerate}
\ep
\begin{proof}
The outline of the proof is as follows. We fix a double structure
$Z$ of type $(a)$ on the line $L$. We have seen that a quasiprimitive
triple line of type $(a,b)$ containing $Z$ is defined via a morphism
$\ideal{Z} \ra \coo_{L}(2a+b)$.  We construct a family of triple lines
$W_t$  specializing  such  a  morphism  to  a  morphism  $\ideal{Z}  \ra
\coo_{L}(-a-2)  \oplus \coo_D$ where  $D$ is  an effective  divisor on
$L$.  The general triple line $W_t$ in the family is
quasiprimitive of type  $(a,b)$, while the special fibre  $W_0$ is the
first  infinitesimal neighborhood  of $L$  plus embedded  points along
$D$.   Finally, we  construct the  desired family  $C_t$ by  picking a
morphism  $\ideal{W_t} \ra  \coo_{L} (3a+c)$  in such  a way  that the
``extra line'' in  $C_0$ covers the embedded points  of $W_0$, so that
$C_0$ is locally Cohen-Macaulay.

We fix coordinates  so that $\pso = \mbox{Proj}  (k[x,y,z,w])$.
Let $L_{0}$ be the line of equations $x=y=0$,
and let $Z_{0}$ be the double structure on $L_{0}$ defined by the
homogeneous ideal $(x^2,xy,y^2,xg-yf) $ where
$f=z^{a+1}$ and  $g=w^{a=1}$. Thus $Z_{0}$ is a  double line of
genus $-a-1$.

Now, over $\Aone = \mbox{Spec}(k[t])$, we consider the
trivial families $L = L_{0} \times \Aone$ and $Z =Z_{0} \times \Aone$.
Let $\ideal{Z}$ denote the ideal sheaf of $Z$ in
$\pso_{\Aone}$. The epimorphism  $\pi:   \ideal{Z}   \rrr
\coo_{L}(-a-2) \oplus \coo_{L} (2a)$  mapping $xg-yf$ to $(1,0)$
and $x^2$, $xy$  and $y^2$ to $(0,f^2)$,$(0,fg)$, $(0,g^2)$
induces an isomorphism $\ideal{Z} \otimes \coo_{L} \cong
\coo_{L}(-a-2) \oplus \coo_{L} (2a)$. The matrix $$M=
\begin{bmatrix} t w^{b}\\ z^{3a+b+2}
\end{bmatrix}
$$ defines a morphism $$\psi: \coo_{L} (-a-b-2) \rrr
\coo_{L}(-a-2) \oplus \coo_{L} (2a),$$ and we  let $\cag = Coker
(\psi)$.  Note that $\cag$ is  flat over $\Aone$ because  $\psi$
remains  injective on  the  fibres over  $\Aone$. Now  let $\chi:
\ideal{Z} \rrr \cag$ be the surjective morphism $$ \ideal{Z}
\stackrel{\pi}{\rrr} \coo_{L}(-a-2) \oplus \coo_{L} (2a) \rrr
\cag $$ and define $W$ by letting $\ideal{W} = Ker
(\chi)$.

It is clear that  $W$ is a flat family of curves  over $\Aone$. For $t
\neq 0$, $\cag_{t}  = \cag \otimes k(t)$ is  isomorphic to $\coo_{L_0}
(2a+b)$,  so that  $W_{t}$ is  a  quasiprimitive triple  line of  type
$(a,b)$.    On   the  other   hand,   $\cag_{0}$   is  isomorphic   to
$\coo_{L_{0}}(-a-2)  \oplus \coo_{D}$,  where  $D$ is  the divisor  in
$L_{0}$ defined by the  global section $s=z^{3a+b+2}$ of $\coo_{L_{0}}
(3a+b+2)$.   By  construction   $W_{0}$  contains  $L_{0}^{(2)}$,  and
$\ideal{L_{0}^{(2)},W_{0}} \cong \coo_{D}$.

We now claim that, over the coordinate ring $R=k[t][x,y,z,w]$ of
$\pso_{\Aone}$, the saturated ideal $I_W$ of $W \subset
\pso_{\Aone}$  has a free graded resolution $$ 0 \rrr F_3
\stackrel{M_{3}}{\rrr} F_2 \stackrel{M_{2}}{\rrr} F_1
\stackrel{M_{1}}{\rrr} F_0 = R $$ where $$
\begin{array}{ccl}
F_1   &=& R(-3)^{\oplus 4} \oplus R(-a-3)^{\oplus 2} \oplus
R(-a-b-2)
\\

\\
F_2   &=& R(-4)^{\oplus 3} \oplus R(-a-4)^{\oplus 4} \oplus
R(-a-b-3)^{\oplus 2}
\\

\\
F_{3} &=& R(-a-5)^{\oplus 2} \oplus R(-a-b-2)
\\
\end{array}
$$ and the maps are defined by the matrices $$
%\begin{array}{ccl}
M_{1}
%&
=
%&
\begin{bmatrix}
x^3 & x^2y  & xy^2 & y^3  & x(xg-yf) & y(xg-yf) &
x^2z^{a+b}+tw^{b} (xg-yf)
\end{bmatrix}
%\\
$$
%\\
$$\begin{array}{cccccc} M_{2} &
=
&
\begin{bmatrix}
y & 0 & 0 &-g & 0 & 0 & 0 & z^{a+b} & 0 \\ -x& y & 0 & f &-g & 0 &
0 & 0 & z^{a+b} \\ 0 &-x & y & 0 & f &-g & 0 & 0 & 0 \\ 0 & 0 &-x
& 0 & 0 & f & 0 & 0 & 0 \\ 0 & 0 & 0 & x & y & 0 &-y & tw^{b} & 0
\\ 0 & 0 & 0 & 0 & 0 & y & x & 0 & tw^{b} \\ 0 & 0 & 0 & 0 & 0 & 0
& 0 & -x & -y
\end{bmatrix}
%\\
&
%\\
M_{3} &
=
&
\begin{bmatrix}
g & 0 &z^{a+b}  \\ -f & g & 0  \\ 0 & -f & 0  \\ y &  0 & 0  \\ -x
& y & 0  \\ 0 & -x & 0  \\ 0 & y & -tw^{b} \\ 0 & 0 & -y \\ 0 & 0
& x
\end{bmatrix}

\end{array}.
$$

We use the Buchsbaum-Eisenbud criterion~\cite{beis-exact} to see that
the complex
above is acyclic. For this it is enough to observe that the ideal
$I_{3} (M_{3})$ generated by the $3 \times 3$ minors of $M_{3}$
contains the regular sequence $(x^3,y^3,z^{3a+b+2})$ and that
$I_{6} (M_{2})$ contains the regular sequence $(x^6,y^6)$.
Therefore the complex is a resolution of the ideal $J$ defined by
$M_{1}$, hence $J$ is the saturated homogeneous ideal of a closed
subscheme of $\pso_{\Aone}$. It is clear that $J \subseteq I_{W}$,
and comparing the Hilbert polynomials we see $J=I_{W}$. This
proves the claim.

We now construct a family $C$ of $4$-lines containing $W$ by letting
$\ideal{C}$ be  the kernel of  a morphism  $\phi: \ideal{W}
\rrr \coo_{L}  (3a+c)$ chosen in such  a way that $C_0$  does not have
embedded points. Recalling that $f=w^{a+1}$ and $g=z^{a+1}$,
the matrix $$
\begin{bmatrix}
        -tf^{3}z^{c-b}w^{b} &
-tf^{2}gz^{c-b}w^{b} & -tfg^{2}z^{c-b}w^{b} & -tg^{3}z^{c-b}w^{b}
& z^{4a+c+3} & z^{3a+c+2}w^{a+1} & w^{4a+b+c+2}
\end{bmatrix}
$$ defines a map $F_{1} \stackrel{P}{\ra} R(3a+c)$ which satisfies
$P M_{2} = 0 \mod I_{L}$, where $I_{L} = (x,y)$ is the ideal of
$L$. We thus obtain a map $I_{W} \stackrel{{\overline P}}{\ra}
T(3a+c)$ where $T = R/I_{L} \cong A[z,w]$. Sheafifying, we obtain
a morphism $\phi: \ideal{W} \ra \coo_{L}(3a+c)$ and we define
our family $C$ by $\ideal{C} = \ker \phi$. Since $\phi_{t}$ is
surjective for all $t \in \Aone$, we see from the proof of
Proposition \ref{qp} above that $C_{t}$ is a quasiprimitive 4-line
of type $(a,b,c)$ for $t \neq 0$.

Finally, we need to check that $C_0$ is a thick 4-line.
To this end, observe that
$$
\ideal{L_0}^{3} \subset \ideal{W_0} \subset \ideal{L_0}^{2}
$$
and the morphism $\phi_{0}$ is zero on
$\ideal{L_0}^{3}$, so  it factors through
$\ideal{W_{0}}/\ideal{L_0}^{3}$. On the other hand, looking at the
presentation   of    $I_{W}   \otimes    T$,   we    see   that
$\ideal{W_{0}}/\ideal{L_0}^{3}$ is isomorphic to
$$\coo_{L_0}(-a-3)^{\oplus 2} \oplus \coo_{L_0}(-a-b-2).$$ Thus we
have a commutative diagram:
\begin{equation}
\begin{CD}
0 @>>> \ideal{C_0}/\ideal{L_0}^{3}  @>>>
\ideal{W_{0}}/\ideal{L_0}^{3} @>\bar{\phi_{0}}>> \coo_{L_{0}}
(3a+c) @>>> 0  \\
      && @VVV @V{\alpha}VV @VVz^{3a+b+2}V \\
0 @>>> \ideal{\bar{C}}  /   \ideal{L_0}^{3}       @>>>
\ideal{L_0}^{2}/\ideal{L_0}^{3} @>{\beta}>> \coo_{L_{0}}
(6a+b+c+2) @>>> 0
\end{CD}
\end{equation}
where, identifying $\ideal{W_{0}}/\ideal{L_0}^{3}$ with
$\coo_{L_0}(-a-3)^{\oplus 2} \oplus \coo_{L_0}(-a-b-2)$ and
$\ideal{L_0}^{2}/\ideal{L_0}^{3}$  with  $\coo_{L_0}(-2)^{\oplus
3}$, the morphisms are $$
\begin{array}{ccl}
\alpha &=&
       \begin{bmatrix} w^{a+1} & 0 & z^{a+b} \\
                          -z^{a+1} &w^{a+1}& 0 \\
               0 & -z^{a+1} & 0
           \end{bmatrix} \\ \\
\beta &=&
       \begin{bmatrix} z^{2a+2} w^{4a+b+c+2} &
           z^{a+1}w^{5a+b+c+3}-z^{6a+b+c+4} &
           w^{6a+b+c+4}-2z^{5a+b+c+3}w^{a+1}
           \end{bmatrix} \\ \\
\bar{\phi_{0}} &=&
       \begin{bmatrix} z^{4a+c+3} & z^{3a+c+2} w^{a+1} & w^{4a+b+c+2}
           \end{bmatrix}
\end{array}
$$ Since $\alpha$ is injective and $$ \deg
(\ideal{L_0}^{2}/\ideal{L_0}^{3}) - \deg (
\ideal{W_{0}}/\ideal{L_0}^{3}) = 3a+b+2, $$ we see that
$C_{0}=\bar{C}$ is a thick 4-line, and this concludes the proof.
\end{proof}

\br\label{qplinespec}{\em It is not known what kind of
specializations might occur between quasiprimitive four-lines,
however we will at least note a necessary condition. If a family
of quasiprimitive four-lines of type $(a,b,c)$ specializes to a
four-line of type $(a^{\prime},b^{\prime},c^{\prime})$, then
$a^{\prime} \leq a$. To see this, consider the deformation of the
underlying double line $Z$: the general such $Z$ has genus $-1-a$,
hence the limit $Z^{\prime}$ consists of a double line of genus
$-1-a^{\prime}$ and possibly some embedded points. Since the
arithmetic genus is constant, we conclude that $a^{\prime} \leq a$.
By the same reasoning, if a family of triple lines of
type $(a,b)$ specializes to another of type $(a^{\prime},
b^{\prime})$, then $a^{\prime} \leq a$. That this actually happens
can be seen in \cite[3.6 and 3.10]{nthree}. \em}\er

\section{A deformation on the double quadric} \label{sdouble}

The goal of this section is to show that families of disjoint
unions of double lines contain certain families of 4-lines in
their closure. This follows readily from our study of curves on
double surfaces \cite{2f}.

Let $F$ be a smooth surface on a smooth projective threefold $T$
and let $X \subset T$ be the effective divisor $2F$. For a curve
$C \subset X$, let $P$ be the curve part of the scheme-theoretic
intersection $C \cap F$. We may write $$ \ideal{C \cap F,F} =
\ideal{Z,F} (-P) $$ where $Z$ is zero-dimensional. The inclusion
$P \subset C \cap F$ generates the commutative diagram
\begin{equation} \label{diagram}\begin{array}{lllllllll}
          & & 0 & & 0 & & 0 & & \\
          & & \downarrow & & \downarrow & & \downarrow & & \\
          0 & \ra & \ideal{R}(-F) & \stackrel{f}{\ra} & \ideal{C} & \ra &
          \ideal{Z,F}(-P) & \ra & 0 \\
          & & \downarrow & & \downarrow & & \downarrow & & \\
          0 & \ra & \coo_{T}(-F) & \stackrel{f}{\ra} & \ideal{P} & \ra &
          \coo_{F}(-P) & \ra & 0 \\
          & & \downarrow & & \downarrow {\phi} & & \downarrow & & \\
          0 & \ra & \coo_{R}(-F) & \stackrel{\sigma}{\ra} & \mathcall &
\ra & \coo_{Z}(-P) & \ra & 0 \\
          & & \downarrow & & \downarrow & & \downarrow & & \\
          & & 0 & & 0 & & 0 & &  \\
\end{array}\end{equation}
which defines the residual curve $R$ to $C$ in $F$. Thus we obtain
a triple $T(C)=\{Z,R,P\}$ in which $R \subseteq P$ are effective
divisors on $F$. Using depth arguments and results on
generalized divisors \cite{hgd} as in \cite{2h,2f}, one finds that
$Z$ is a Gorenstein divisor on $R$, $\mathcall \cong \coo_{R}(Z-F)$
is a rank one
reflexive $\coo_{R}$-module, and $\sigma$ gives a section of
$\mathcall(F)$ that defines $Z$ as a generalized divisor on $R$.
Note that the arithmetic genus of $C$ is given by the
formula
\begin{equation}
g(C)= g(P) + g(R) + \deg_R \coo_{R}(F) - \deg (Z) -1.
\end{equation}

In the following lemma we consider the case $X$ is  a double quadric
surface $2Q$ in $\pso$, and describe the triple $T(C)$ for
a general quasiprimitive $4$-line $C$ of type $(0,b,c)$. This will
allow us to conclude that $C$ is a specialization of a disjoint union
of double lines, a fact we  will later use
in our description of the
irreducible components of the Hilbert scheme (Theorem \ref{comps})
and in the proof of connectedness of the Hilbert scheme
(Theorem~\ref{hilbconnected}).

\bl \label{sticky} Let $C \subset \pso$ be a general
quasiprimitive $4$-line of type $(0,b,c)$ with Cohen-Macaulay
filtration  $L \subset D \subset W \subset C$. Then there exists
a smooth quadric surface $Q$ such that
\begin{enumerate}
       \item $D \subset Q$ and $C \subset 2Q$.
       \item The triple associated to $C \subset 2Q$ has form $T(C) = 
\{Z,D,D\}$,
       where $Z$ consists of $c-b$ simple points and $b+2$ double points,
       none of which are contained in $L$.
\end{enumerate}
Further, $H^{1}(\coo_{D}(Z+D-Q))=0$.
\el

\begin{proof} Existence of the smooth quadric in the first statement
follows from \cite[1.5]{nthree}. Fixing homogeneous coordinates
$x,y,z,w$ on $\pso$ so that $x$ and $y$ generate the ideal of $L$,
we may assume that $F=xw-yz$ is an equation for $Q$. Then $D$ is
the divisor $2L$ on $Q$, and the isomorphism $\ideal{L,D} \cong
\coo_L$  maps the global sections of $\ideal{L,D}(1)$ defined by
$x$ and $y$ to $z$ and $w$ respectively, where we consider $z$ and
$w$ as elements of $\text{H}^{0} (L, \coo_L (1))$.

As explained in Proposition~\ref{qp} and \cite[2.3]{nthree}, $W$
corresponds to an epimorphism $$ \psi: \ideal{D}/ {\ideal{L}
\ideal{D}} \cong \coo_L  \oplus \coo_L(-2) \stackrel{[\, p \;\;q
\,]}{\longrightarrow} \coo_{L} (b). $$ in which $\psi$ maps the
global section defined by $F$ to $q \in \text{H}^0 (L, \coo_L
(b+2))$. Thus under the isomorphism $\ideal{D,W} \cong \coo_L
(b)$, the section defined by $F$ is mapped to $q$. Since $W$ is
not contained in $Q$ by assumption, we see that $q$ is nonzero.

Clearly $W$ is contained in $2Q$ and its associated triple has the
form $(Z_W,L,D)$ for some zero dimensional subscheme $Z_W \subset
L$ of degree $b+2$. In fact, if $q$ is an equation for $Z_W$, then
the bottom row of diagram~\ref{diagram} becomes $$ 0 \ra \coo_{L}
(-2) \stackrel{q}{\ra} \coo_L (b)\cong \ideal{D,W} \ra \coo_{Z_W}
(-D) \ra 0 $$ with the identifications above. Putting this
together with the analogous sequence for the triple $\{Z,D,D\}$
associated to $C$ we obtain the commutative diagram:
       $$\begin{array}{lllllllll}
& & 0 & & 0 & & 0 & & \\ & & \downarrow & & \downarrow & &
\downarrow & & \\ 0 & \ra & \ideal{L,D} (-2) \cong \coo_L (-2) &
\stackrel{\alpha}{\ra} & \ideal{W,C} \cong \coo_L(c) & \ra &
\ideal{Z_W,Z}(-D) & \ra & 0 \\ & & \downarrow & & \downarrow &  &
\downarrow & & \\ 0 & \ra & \coo_{D} (-2) & \ra & \ideal{D,C} &
\ra & \coo_{Z_{\;}}(-D) & \ra & 0 \\
            & & \downarrow & & \downarrow & & \downarrow & & \\
0 & \ra &  \coo_{L} (-2) & \ra &
            \ideal{D,W} \cong \coo_L (b)  & \ra & \coo_{Z_W} (-D) & \ra & 0 \\
            & & \downarrow & & \downarrow & & \downarrow & & \\
            & & 0 & & 0 & & 0 & &  \\
       \end{array}$$
We need to describe the first horizontal arrow $\alpha$ in the
diagram. First note that $\alpha$ factors through the inclusion
$\ideal{L,D} (-2) \cong \coo_L (-2) \stackrel{\beta}{\ra}
{\ideal{L} \ideal{D}}/ \mathcal{J} \cong \coo_{L} (b)$ where
$\mathcal{J} = \ideal{L}\ideal{W} + \ideal{Z}^2$ as in
Proposition~\ref{qp}. The map $\beta$ in turn is obtained from the
above inclusion $\coo_L (-2) \stackrel{q}{\ra} \coo_L (b)\cong
\ideal{D,W}$ upon tensoring with $\ideal{L,D} \cong \coo_L$.
Finally, ${\ideal{L} \ideal{D}}/ \mathcal{J} $ is the $\coo_L (b)$
summand in the decomposition $\ideal{W}/\mathcal{J} \cong \coo_{L}
(b) \oplus \coo_{L} (-b-2)$ of Proposition~\ref{qp}, and the
isomorphism $\ideal{W,C} \cong \coo_L (c)$ arises from the exact
sequence $$ 0 \ra \ideal{C}/\mathcal{J} \ra \ideal{W}/\mathcal{J}
\cong \coo_{L} (b) \oplus \coo_{L} (-b-2) \stackrel{[\, f \;\; g
\,]}{\longrightarrow} \coo_{L} (c) \ra 0. $$ We conclude that
$\alpha$ is multiplication by the global section $fq \in
\text{H}^{0} (L, \coo_L (c+2))$, and therefore $\ideal{Z_W,Z}
\cong \coo_L / (fq)$.

In the irreducible family of quasiprimitive $4$-lines of type
$(0,b,c)$, the set of curves $C$ for which the form $qf$ has
simple zeros is open. Looking at $$ \exact{\ideal{Z_W,Z} \cong
\coo_L / (fq)}{\coo_Z }{\coo_{Z_W} \cong \coo_L / (q)} $$ we see
that subscheme $Z$ associated to a general such curve $C$ consists
of one double point for each zero of $q$, and one simple point for
each zero of $f$; Moreover, $qf=0$ defines the divisor $Z \cap L
\subset L$.

Finally, we show that $H^{1}(\coo_{D}(Z+D-Q))=0$.
Consider $W = Z \cap L $ and the residual
scheme $Y$ to $W$ in $Z$ From the description above, $\deg W =
c+2$ and $\deg Y = b+2$. Since $D^{2}=0$ on $Q$, the sequence
relating $Y$ and $W$ to $Z$ takes the form $$ 0 \rrr \ideal{Y,L}
\rrr \ideal{Z,D} \rrr \ideal{W,L} \rrr 0 $$ and applying
$\Hom_{\coo_{D}}(-,\coo_{D})$ yields the exact sequence $$ 0 \rrr
\coo_{L}(W) \rrr \coo_{D}(Z) \rrr \coo_{L}(Y) \rrr 0 $$ (one
checks locally that $\Ext^{1}_{\coo_{D}}(\coo_{L},\coo_{D})=0$ and
$\Hom_{\coo_{D}}(\coo_{L}(a),\coo_{D}) \cong \coo_{L}(-a)$ by
\cite[III,6.7]{AG}). Tensoring by $\coo_{D}(-Q)$ and taking the
long exact cohomology sequence now gives the desired vanishing.
\end{proof}

\bp \label{disjointdoublelines} Let $C$ be a general quasiprimitive
4-line of type $(0,b,c)$ lying on a double quadric surface $2Q$
as in Lemma~\ref{sticky}. Then there is a flat family of
curves $X_{t}$, $t \in \Pone$ on $2Q$ with $X_{0}=C$ and general
member a disjoint union of double lines of arithmetic genera
$-1-b$ and $-1-c$. \ep

\begin{proof}
Since the double line $D$ underlying $C$ is linearly
equivalent on $Q$ to the disjoint union of two lines, our
claim is a consequence of Lemma~\ref{sticky} and of the
following theorem from~\cite{2f}
\end{proof}

\bt\label{limits} Let $C \subset X$ be a curve with triple
$T(C)=\{Z,R,P\}$ such that $H^{1}(\coo_{R}(Z+P-F))=0$. Then $C$ is
a specialization of curves $C^{\prime}$ with triples
$\{Z^{\prime},R^{\prime},P^{\prime}\}$ for which $R^{\prime}$ is
general in its divisor class on $F$. \et

\section{Curves of degree $4$ with large speciality} \label{twisted}
In this section we study curves which have large speciality. We
express the speciality of a curve $C$ by its {\it spectrum}
which can be defined~\cite{sch-tran} as the non-negative function
$$
h_{C}(n)=\Delta^{2} h^{0}(\coo_{C}(n))=
h^{0}(\coo_{C}(n)) - 2 h^{0}(\coo_{C}(n-1)) + h^{0}(\coo_{C}(n-2)),
$$
which we represent by the t-uple
of integers with exponents $\{n^{h_{C}(n)}\}$.
The curves with the largest speciality are the extremal curves, which form
an irreducible component $E \subset H_{4,g}$ of dimension $15-2g$
\cite[Theorem 2.5 and Theorem 3.7]{extremal}. Extremal curves are defined
as those curves achieving upper bounds \cite{bounds} on the Rao function
$h^{1} (\ideal{C} (n))$,
but for $d=4$ and $g \leq 0$ they may also be characterized as
(a) nonplanar curves containing a plane cubic curve, or
(b) curves with spectrum $\{g\} \cup \{0,1,2\}$ \cite[2.2]{subextremal}.

Similarly, there are sharp upper bounds on the Rao function for
non-extremal curves~\cite{subextremal}. The curves achieving these
bounds are called
subextremal and have spectrum $\{g+1,0,1^{2}\}$, although they are
not characterized by this fact \cite[2.15]{subextremal}. The curves
with the speciality of a subextremal curve are characterized as follows.

\bl \label{character}Let $C$ be a curve of degree $4$ and genus $g
\leq -2$. Then $C$ has spectrum $\{g+1,0,1^{2} \}$ if and only if
$C$ contains a subcurve $T$ of degree $3$ and genus $0$.
\el

\begin{proof}
If $C$ contains a curve $T \in H_{3,0}$, then the principal spectrum
spectrum $\{0,1^{2}\}$ of $T$ is contained in that of $C$
\cite[\S 3]{sch-tran} and the remaining element $g+1$ is determined
by the genus of $C$.

Conversely, suppose that $C$ has spectrum $\{g+1,0,1^{2} \}$. Then
$h^0 \coo_{C}(g+1) =1$ and choosing $0 \neq \alpha \in H^0 \coo_{C}(g+1)$
gives a map $\coo_{C} \rrr \coo_{C}(g+1)$
with image $\coo_{D}$ for some closed subscheme
$D \subset C$. The local depth of $D$ is one
because $\coo_{D} \subset \coo_{C}(g+1)$, hence $D$ is a locally
Cohen-Macaulay curve.
The inclusion above also shows that $h^0 \coo_{D} (1) \leq
2$: for $g < -2$ this is because $h^{0} \coo_{C}(g+2)=2$. If $g=-2$,
then $h^0 \coo_{C} = 3$ and the inclusion
$\mbox{H}^{0} (\coo_{D} (1)) \subset \mbox{H}^{0} (\coo_{C})$ is
strict, as otherwise we obtain a surjection $\coo_{D}(1) \ra \coo_C$,
which is absurd. It follows that $D=L$ is a line.
This yields an exact sequence

\begin{equation}\label{J}
0 \rrr \coo_{L}(-g-1) \rrr \coo_{C} \rrr \coo_{T} \rrr 0
\end{equation}
for a closed subscheme $T \subset C$ of degree $3$ and genus $0$.

If $T$ is not purely one-dimensional, then the purely one-dimensional part
$P \subset T$ is planar because $g(P) > 0$ and $\deg P =3$ \cite[3.1]{genus},
but this is not possible because $C$ is not extremal. Thus $T$ is
locally Cohen-Macaulay,
and this finishes the proof.
\end{proof}

\bp \label{coh}
Let $C$ be a curve of degree $4$ and genus $g \leq -1$ having
spectrum $\{-g+1,0,1^{2}\}$.
Then there is a line $L$ such that the Rao module
$M_{C} = H^{1}_{*}(\ideal{C})$ is a graded module over the coordinate ring
$S_{L}$ with resolution
\begin{equation}\label{rm}
0 \rrr S_{L}(-j) \oplus S_{L}(j-5+g) \stackrel{\sigma}{\rrr}
S_{L}(-2)^{\oplus 3} \rrr S_{L}(-g-1) \rrr M_{C} \rrr 0
\end{equation}
for some integer $2 \leq j \leq n(g) = \lfloor \frac{5-g}{2} \rfloor$.
The cohomology of $C$ is determined by $j$ and we denote the
corresponding family of curves by $H_{j}$.
\ep

\begin{proof}
We first treat the case $g \leq -2$. In this case $C$ contains
a degree $3$ and genus $0$ curve $T$ by~\ref{character}.
As every curve of degree $3$ and genus $0$,  $T$ is arithmetically
Cohen-Macaulay with total ideal $I_T$ generated by three
quadrics \cite[3.5]{genus}. Furthermore, by the proof of~\ref{character},
there is a line $L$ such that $\ideal{T,C} \cong \coo_{L} (-g-1)$.
Factoring the surjection $\coo_{\pthree}(-2)^{\oplus 3} \rrr \ideal{T}
\rrr \coo_L(-g-1)$ through $\coo_L(-2)^{\oplus 3}$ and writing the
kernel of the induced map as $\coo_{L}(-j) \oplus \coo_{L}(j-5+g)$
for some integer $j \in \lfloor 2,\frac{5-g}{2} \rfloor$, we
obtain the resolution (\ref{rm}).

Now assume that $g=-1$. The spectrum shows that $C$ is
neither ACM (because $h_{C}(0) > 1$) nor extremal, hence
$h^{1} \ideal{C} \leq 1$, $h^{1} \ideal{C}(1) \leq 2$ and
$h^{1} \ideal{C}(2) \leq 1$ \cite[Theorem 2.11]{subextremal}.
The first two of these are equalities in
view of the Euler characteristics (since $h^{3} \ideal{C} = h^{3}
\ideal{C}(1)=0$) and $h^{2} \ideal{C}(1)=0$.
If $h^{1} \ideal{C}(2)=1$, then $C$ is subextremal by definition and
resolution (\ref{rm}) for $j=2$.

If $h^{1} \ideal{C}(2)=0$, then $\ideal{C}$ is $3$-regular
(hence $h^{1} \ideal{C}(n)=0$ for $n \geq 2$)
and $h^{0} \ideal{C}(2)=h^{2} \ideal{C}(2)=0$.
Since $h^{1} \ideal{C}(n)$ is increasing or zero for $n \leq 0$, we
see that $h^{1} \ideal{C}(n)=0$ for $n < 0$.
In particular, the Rao module $M_{C}$ has a generator $m$ in degree $0$.
If $m$ does not generate $M_{C}$ as a
module over the homogeneous coordinate ring $S=k[x,y,z,w]$ of $\pso$,
then $m$ is annihilated by $3$ independent linear forms,
which implies that $C$ lies on a quadric by \cite[3.4.5]{sch}, a
contradiction. Thus $m$ generates the Rao module and
$M_{C} \cong S/(x,y,z^{2},zw,w^{2})$ after a change of coordinates, so
that $M_{C}$ has resolution (\ref{rm}) for $I_{L}=(x,y)$ and $j=3$.
\end{proof}

Finally, we describe the families $H_{j}$ and how they fit
together in the Hilbert scheme.

\bp \label{d3g0-c} For fixed $g \leq -1$, let $H_{j}$ denote the
family of curves $C \in H_{4,g}$ with spectrum $\{g+1,0,1^{2}\}$ and
Rao module $M_{C}$ having resolution (\ref{rm}).
Let $G_4$ denote the family of thick 4-lines. Then
\begin{enumerate}
\item
The family $H_{2}$ is irreducible of dimension $13-2g$ (resp. $16$ if $g=-1$).
It is the family of subextremal curves and meets $G_4$.
\item
The family $H_{3}$ is irreducible of dimension $13-2g$ (resp. $16$ if $g=-1$).
%The general curve in $H_{3}$ is the union of a
%smooth conic and a double line of genus $g-1$ meeting in a scheme of
%length two.
It meets $G_4$ in general and $G_4 \subset {\overline {H_3}}$ if $g \geq -2$.
\item
Suppose $g \leq -3$. Then $H_{j} \subset G_4$ for $3 < j \leq n(g)$, and
$G_4 = {\overline {H_{n(g)}}}$ is irreducible of dimension $9-3g$.
\end{enumerate}
\ep

\begin{proof}
We consider the last statement first: suppose $g \leq -3$ and let
$C \in H_{j}$ for $3<j\leq n(g)$.
The sequence (\ref{rm}) shows that $h^{1} \ideal{C}(3)=-g-3$, hence
$h^{0} \ideal{C}(3)=4$. If $L$ and $T$ are as in the proof of
\ref{character} above, then $\dim (I_{L} I_{T})_{3} \leq 4$
because $I_{L} I_{T} \subset I_{C}$ and so $I_{T}=I_{L}^{2}$ by
Lemma \ref{powers} below. It follows that $C$ is a thick 4-line
supported on $L$, hence $H_{j} \subset G_4$.
As we saw in Proposition \ref{thickmoduli}, $G_4$ is irreducible and
the thick 4-lines supported on $L$ are parametrized by the open subset
$$U \subset {\rm Hom}_{\coo_{\pthree}} (\ideal{T},\coo_{L}(-g-1))
\cong {\rm Hom}_{\coo_{L}} (\coo_{L}(-2)^{3}, \coo_{L}(-g-1))$$
corresponding to surjective maps. For $j=2$ and $3$, we can use the specific
surjection given by $(w^{1-g},w^{3-g-j}z^{j-2},z^{1-g})$ to see that
$G_4$ meets $H_{j}$. Since $I_{C}$ is determined by its image
in $I_{T}/{I_{L} I_{T}} \cong S_{L}(-2)^{3}$ and hence by the image
of the first map in sequence (\ref{rm}) above, we find by counting
dimensions that $H_{j}$ is irreducible of dimension $5+2j-2g$
(except if $2j=5-g$, when the dimension is $4+2j-2g$). For
$j=n(g)=\lfloor \frac{5-g}{2} \rfloor$, the closure of $H_{n(g)}$
is irreducible of dimension $9-3g$, hence is equal to $G_4$.

If $j=2$, then $H_{j}$ consists of subextremal curves, which are
those obtained from extremal curves of degree $2$ and genus
$g^{\prime}=g-1$ by a height one biliaison on a quadric surface
\cite[2.11 and 2.14]{subextremal}.
Let then $\gamma$, $\rho$ (resp. $\gamma^{\prime}$, $\rho^{\prime}$) be the
gamma and Rao functions for the extremal curves of degree $2$ and
genus $g-1$ (resp. subextremal curves of degree $4$ and genus
$g$). Letting $\B_{\gamma,\rho,2,1}$ denote the universal biliaison scheme
of Martin-Deschamps and Perrin \cite[VII \S 4]{MDP}, we have smooth
irreducible projections
$$
\begin{array}{ccc}
        \B_{\gamma,\rho,2,1} & \stackrel{q_{2}}{\ra} &
H_{\gamma^{\prime},\rho^{\prime}}  \\
        q_{1} \downarrow   & & \\
E=H_{\gamma,\rho} & &
\end{array}
$$
to the spaces $H_{\gamma,\rho}$ (resp. $H_{\gamma^{\prime},\rho^{\prime}}$)
of curves with constant cohomology. The family
$E=H_{\gamma^{\prime},\rho^{\prime}}$ of extremal curves is irreducible
of dimension $7-2g$ \cite[2.5]{extremal} and using \cite[VII, 4.8]{MDP}
we compute that the fibre dimension of $q_{1}$ is $8$ (resp. $9$ if
$g=-1$) and the fibre dimension of $q_{2}$ is $2$, hence
the family $H_{\gamma,\rho}$ of subextremal curves is
irreducible of dimension $13-2g$ (resp. $16$ if $g=-1$).

For $j=3$, we take an indirect approach.
Consider the family of arithmetically Cohen-Macaulay curves $D$ with resolution
of the form
$$0 \ra \coo(2g-1) \oplus \coo(g-2) \ra \coo(g-1)^{3} \ra \ideal{D}
\ra 0.$$
This family is irreducible of dimension
$2 \binom{-g+3}{3}+ \binom{-g+2}{2}+2$ (resp. $12$ if $g=-1$) by \cite{Ell}
and the general member is smooth and irreducible
(the numerical character has no gaps).

Let $D$ be a general such curve. Then a general (disjoint) union $D \cup L$
with a line $L$ lies on an integral surface of degree $-g+2$; To see this,
use the linear systems $\PP H^{0} \ideal{D} (-g+1)$ and
$\PP H^{0} \ideal{L} (1)$ to obtain a map $\tau: \pthree-D-L \ra \ptwo \times
\pone$. Since the fibres of $\tau$ (= $S \cap H:S \in H^{0} \ideal{D} (-g+1),
H \in H^{0} \ideal{L} (1)$) are generally of dimension one,
the image of $\tau$ has dimension two.
Composing with the Segre embedding $\ptwo \times \pone
\hookrightarrow {\PP}^{5}$, we apply Jouanolou's Bertini theorem
\cite[6.10]{J} to see that the general surface of degree $-g+2$
containing $D \cup L$ is irreducible. Furthermore, the resolution for
$\ideal{D}$ shows that
$H^{1} \coo_{D \cup L} ((-g+2)-4-(g+1)) \neq 0$
and we find that \cite[III, 2.7(b)]{MDP} $D \cup L$ can be
bilinked on a surface of degree $-g+2$ with height $g+1$ to a curve
$C$, which lies in $H_{3}$ by direct calculation.

Let $\gamma, \rho$ (resp. $\gamma^{\prime}, \rho^{\prime}$) be
the gamma and Rao functions for curves in $H_{3}$ (resp. $D \cup L$).
Letting $\B_{\gamma,\rho,-g+2,-g-1}$ be the universal biliaison scheme
\cite[VII \S 4]{MDP}, we obtain smooth irreducible projections
$$
\begin{array}{ccc}
        \B_{\gamma,\rho,-g+2,-g-1} & \stackrel{q_{2}}{\ra} &
H_{\gamma^{\prime},\rho^{\prime}}  \\
        q_{1} \downarrow   & & \\
H_{3}=H_{\gamma,\rho} & &
\end{array}$$
   From the last paragraph, the image of $q_{2}$ is dense in the irreducible
component consisting of the closure of the family of
disjoint unions $D \cup L$ considered above. Using the resolutions
given, we compute the dimension of the fibres of $q_{1}$ and $q_{2}$
via \cite[VII 4.8]{MDP} and conclude that $H_{3}$ is
irreducible of dimension $13-2g$ (resp. $16$ if $g=-1$).
\end{proof}

\br\label{descriptions}{\em
One can check by a dimension count that the
general members of the families $H_{2}$ and $H_{3}$ are described
as follows.
\begin{enumerate}
       \item For $g=-1$, the general member of $H_{2}$ is a disjoint
       union of conics. For $g \leq -2$, the general member of $H_{3}$
       is the union of  a double line $Z$ of genus $g-2$ and two disjoint
       lines $L_1$ and $L_2$, each meeting $Z$ in a scheme of length $2$.

       \item For $g=-1$, the general member of $H_{3}$ is a disjoint
       union of a line and a twisted cubic curve. For $g \leq -2$, the
       general member of $H_{3}$ is the union of a double line  of genus $g-1$
       and a smooth conic  meeting in a scheme of length $2$.

\end{enumerate}
\em}\er

The following lemma and its proof are well known:
\bl\label{powers}
Let $L \subset \psn$ be a linear subvariety of codimension two and let
$I$ be the homogeneous ideal generated by a subspace
$V \subset \mbox{H}^{0} (\psn, \coo (d))$ of dimension $r$.
Then the image $W$ under the multiplication map
$V \otimes \mbox{H}^{0}(\psn, \ideal{L} (1)) \rightarrow
\mbox{H}^{0} (\psn, \coo (d+1))$
satsifies $\dim (W) \geq r+1$ with equality if and only $I =I_{L}^{r-1}f$
for some form $f$ of degree $d-r+1$.
\el
%looks like Gotzmann's interpretation of Macaulay's growth bound

\begin{proof}
Let $S={\rm Sym} \; \mbox{H}^{0}(\psn, \ideal{L} (1)) \cong k[x,y]$
be the symmetric algebra and set $\pone = \mbox{Proj} (S)$.
Sheafifying  the natural  map
$V \otimes_{k} S \rrr \bigoplus_{n} \mbox{H}^{0} (\psn,\coo(n))$
of free graded $S$-modules over $\pone$ and letting $\caf$ denote the
image, we obtain an exact sequence
$ \exact{\cae}{\coo_{\pone}^{r}}{\caf}$
of locally free $\coo_{\pone}$-modules. By hypothesis, $h^0(\cae) = 0$ and
$\rank (\cae) \leq r-1$, so we may write
$$\displaystyle \cae \cong \bigoplus_{i=1}^{s} \coo_{\pone} (-a_i)$$
with $s \leq r-1$ and $a_{i} > 0$, hence
$\mbox{h}^{0} (\cae(1)) \leq r-1$.
Since $\mbox{H}^{0} (\cae(1))$ is the kernel of the surjection
$V \otimes S_1 \ra W$,
we see that $\dim (W) \geq r+1$ with equality if and only if
$s=r-1$ and $a_{i} =1$ for $1 \leq i \leq r-1$, which is equivalent
to saying that $\caf \cong \coo_{\pone}(r-1)$.
\end{proof}

\section{Triple lines union a line} \label{striple}
In this section we are interested in families of curves $C$ that are
unions of a quasi-primitive triple line $W$ of type $(a,b)$ and
a reduced line $L$, when $W \cap L$ is non-empty.  Note that
we have
$$g(C) = g(W) + g(L) + \Length  (W \cap L) -1 = -3a -3 -b(W) +
\Length  (W \cap L).$$
In what follows, we will fix $g$ and  $a$ with $a \geq 0$.
If $Z \subset W$ denotes the underlying double line of genus $-1-a$,
we have four different families of such curves $C$ in $H_{4,g}$:
\begin{itemize}
       \item[] $F_{1}=\{W \cup L: \Length    (W \cap L)=3,\,b(W) =-3a-g \}$
       \item[] $F_{2}=\{W \cup L: \Length    (W \cap L)=2$,
       $\Length (Z \cap L)=2, \, b(W) = -3a-g-1 \}$
       \item[] $F_{3}=\{W \cup L: \Length    (W \cap L)=2$,
       $\Length (Z \cap L)=1, \, b(W) -3a-g-1 \}$
       \item[] $F_{4}=\{W \cup L: \Length    (W \cap L)=1, \, b(W) =-3a-g-2\}$
\end{itemize}
The main results of this section are that families $F_{1}$ and $F_{2}$
are irreducible (Proposition \ref{WL}) and that the other two families lie in
their closures (Proposition \ref{WLclosure}).

Our arguments hinge on the following correspondence: Fix a double line
$Z$ of type $a \geq 0$ with support $Y$ and let $b \geq a$ be a integer.
Let $S \supset Z$ be a surface of degree $a+b+2$ with equation $h$
which does not contain the first infinitesimal neighborhood $Y^{(2)}$ of $Y$.
Removing the possible embedded points from the scheme cut out by
$I=(I_{Y}I_{Z},h)$ yields a quasi-primitive $3$-line $W \supset Z$.
This defines $\Phi$:

$$
\mathcals=\{h \in H^{0}(\ideal{Z}(a+b+2)): h \notin I_{Y}^{2}\}
\stackrel{\Phi}{\rightarrow} \{\mbox{Quasiprimitive 3-lines } W \supset Z\}
$$

\bl \label{corresp} Let $\Phi$ be the map above. Then
\begin{enumerate}
\item The image of $\Phi$ is the set of quasiprimitive 3-lines
$W \supset Z$ of type $(a,b^{\prime})$ with $b^{\prime} \leq b$.
\item For $h$ as above, set $I = (I_{Y}I_{Z},h)$. Then the following are
equivalent.
\begin{enumerate}
       \item $W=\Phi(h)$ has type $(a,b)$.
       \item The scheme defined by $I$ is locally Cohen-Macaulay.
       \item $h$ is irreducible modulo $I_{Y}I_{Z}$.
\end{enumerate}
If any of these conditions hold, then $I=I_{W}$.
\end{enumerate}
\el

\begin{proof}
For $h \in (I_{Z})_{a+b+2}: h \notin I_{Y}^{2}$, let $W$ be the
purely one-dimensional part of the scheme $V$ defined by the ideal
$I=(I_{Y}I_{Z},h)$. Since $W$ is quasi-primitive and contains $Z$,
$W$ has type $(a,b^{\prime})$ for some $b^{\prime} \geq a$ and the
total ideal may be written $I_{W}=(I_{Y}I_{Z},h^{\prime})$ with
$h^{\prime}$ of degree $a+b^{\prime}+2$ by \cite[2.3]{nthree}.
The inclusions $(I_{Y}I_{Z},h) \subset I_{V} \subset I_{W}$ now show that
$b^{\prime} \leq b$. On the other hand, if
$W$ is a quasi-primitive 3-line of type $(a,b^{\prime})$ with
$b^{\prime} \leq b$, then
writing $I_{W}=(I_{Y}I_{Z},h^{\prime})$ as above and choosing
a hypersurface $F$ of degree $b-b^{\prime}$ with equation $f$
meeting $Z$ properly, we see that $\Phi(fh^{\prime})=W$.

For the equivalences in statement 2, let $C$ be the scheme defined by
$I=(I_{Y}I_{Z},h)$ so that $W=\Phi(h)$ is obtained from
$C$ by removing possible embedded points.\\
$(a) \Rightarrow (b)$ If $W$ has type $(a,b)$, then by
\cite[2.3]{nthree} the total ideal for $W$ takes the form
$I_{W}=(I_{Y}I_{Z},h^{\prime})$ with $\deg h^{\prime} = a+b+2=\deg h$.
The inclusions $I \subset I_{C} \subset I_{W}$ show that all three ideals
are equal, so $C=W$ is locally Cohen-Macaulay. \\
$(b) \Rightarrow (c)$ Suppose that $h=h^{\prime} t$ modulo $I_{Y}I_{Z}$.
If $T$ is the surface with equation $t$, then we may assume $T$ meets
$Y$ properly (since if both $t \in I_{Y}$ and $h^{\prime} \in I_{Y}$,
then $h \in I_{Y}^{2}$, contrary to hypothesis). In this case
$(I_{Y}I_{Z}, t)$ defines a scheme of length $4 \deg T$
(because $I_{Y} I_{Z}$ defines a locally Cohen-Macaulay 4-line) and
if $C$ is locally Cohen-Macaulay, then $(I,t)$ defines a
scheme of length $3 \deg T$; since $(I_{Y}I_{Z}, t)=(I,t)$, we
must conclude that $C$ is not locally Cohen-Macaulay. \\
$(c) \Rightarrow (a)$ If $W$ has type $(a,b^{\prime})$ for $b^{\prime}<b$,
then $I_{W}=(I_{Y}I_{Z},h^{\prime})$ with $\deg h^{\prime} = a+b^{\prime}+2$.
The inclusion $I \subset I_{W}$ shows that there exists $t$ of degree
$b-b^{\prime}$ such that $h=h^{\prime} t$ modulo $I_{Y}I_{Z}$.

\end{proof}

\bp \label{WL} With the notation above, in $H_{4,g}$ we have
\begin{enumerate}
       \item $F_{1}$ is irreducible of dimension $11-2g-a$
               if $\ds 0 \leq a \leq -\frac{g}{3}$ and
               empty if $\ds a > -\frac{g}{3}$.
       \item $F_{3}$ is irreducible of dimension $10-2g-a$
               if $\ds 0 \leq a < \frac{-g-1}{3}$ and
               empty if $\ds a \geq \frac{-g-1}{3}$.
\end{enumerate}
\ep

\begin{proof}
We first prove statement 1, then indicate the changes to obtain
statement 2.
Let $H \subset H_{3,-a}$ be the family of unions $Z \cup_{2P} L$.
By \cite[3.2(a)]{nthree}, $H$ is irreducible of dimension $9+2a$.
For $b=-3a-g$ we interpret $H^{0} \coo_{\pso}(a+b+2)$ as an affine scheme.
Pulling back the universal family over $H$ we obtain a diagram
$$
\begin{array}{ccc}
\mathcalz \cup \mathcall & \subset & \PP^{3} \times
H \times H^{0}(\coo_{\pthree}(a+b+2)) \\
    & \searrow & \downarrow  \\
& & H \times H^{0}(\coo_{\pthree}(a+b+2)). \\
\end{array}
$$

Consider the closed subset
\begin{equation}\label{v}
       V = \{(Z \cup L,h) \in H \times H^{0}(\coo_{\pthree}(a+b+2)):
h \in I_{Z} \cap (I_{Y}^{3},I_{L})\}
\end{equation}
with first projection $V \stackrel{p_{1}}{\rightarrow} H$.
The fibres of $p_{1}$ are vector subspaces of dimension
$\binom{a+b+5}{3}-3a-2b-7$.
Indeed, after a change of coordinates we may write
$I_{Z}=(x^{2},xy,y^{2},xg-yf)$ (\cite[1.4(c)]{nthree}) and $I_{L}=(x,z)$,
when the fibre is identified with the kernel $K$ of the composite map
\begin{equation}\label{mp}
       H^{0}(\ideal{Z}(a+b+2)) \hookrightarrow (x,z,y^{2})_{a+b+2} \rightarrow
((x,z,y^{2})/(x,z,y^{3}))_{a+b+2} \cong k
\end{equation}
(the inclusion has the correct target because $Z$ meets
$L$ in the double point $2P$). Since $\ideal{Z}$ is $(a+2)$-regular,
$h^{0}(\ideal{Z}(a+b+2))=\chi \ideal{Z}(a+b+2)= \binom{a+b+5}{3}-3a-2b-6$
and $p_{1}:V \ra H$ is an affine bundle with fibres of dimension
$\binom{a+b+5}{3}-3a-2b-7$. In particular, $V$ is irreducible.

Consider the open subset $U=\{(Z \cup L,h) \in V: h \not\in I_{Y}^{2}\}$.
The correspondence of Lemma \ref{corresp} shows that elements of $U$
determine unions $W \cup L$. If $\mathcals \subset \pso \times U$ is
the family of surfaces with equation $h$, we obtain flat families
$$
\begin{array}{ccc}
\mathcaly, \mathcalz \cup \mathcall, \mathcals &
\subset & \pso \times U \\
& \searrow & \downarrow  \\
& & U \\
\end{array}
$$
where $\mathcaly$ is the support of $\mathcalz$.
The subscheme $\mathcalw \subset \PP^{3} \times U$ defined by the ideal sheaf
$\ideal{\mathcalw} = \ideal{\mathcaly} \ideal{\mathcalz} + \ideal{\mathcals}$
is also flat over $U$. To see this, observe that the fibre of the
sheaf $\caf = \ideal{\mathcalw} = \ideal{\mathcaly} \ideal{\mathcalz} +
\ideal{\mathcals}/{\ideal{\mathcaly} \ideal{\mathcalz}}$ over
$u=(Z \cup L,h)$ is isomorphic to $\coo_{Y}(-a-b-2)$: Indeed,
$\caf_{u}$ is a quotient of $\coo_{Y}(-a-b-2)$ via the generator $h$,
and is a torsion free $\coo_{Y}$-module because of the inclusion
$\caf_{u} \subset \ideal{Z}/{\ideal{Y}\ideal{Z}}$. It follows that the
fibres of $\mathcalw$ have constant Hilbert polynomial and the family
is flat by \cite[III, 9.9]{AG}.

Finally, let $U^{\prime} \subset U$ be the open set for which
the fibres of $\mathcalw$ are locally Cohen-Macaulay, taken with the
induced reduced scheme structure.
This is precisely the set for which the fibres of $\mathcalw$ are
quasiprimitive 3-lines of type $(a,b)$ by Lemma \ref{corresp}.
The curves used in the proof of Corollary \ref{extend} show that
$U^{\prime}$ is non-empty; the double line $Z$ with total ideal
$I_{Z}=(x^{2},xy,y^{2},xz^{a+1}-yw^{a+1})$ has the
line $L = \{x=w=0\}$ as a double tangent and
$h=z^{b}(xz^{a+1}-yw^{a+1})-x^{2}w^{a+b}=0$ satisfies the conditions above.
The definition of $V$ above makes it clear that the fibres of
$\mathcalw$ meet the lines $L$ in triple points, so the family
$\mathcalw \cup \mathcall$ is also flat over $U^{\prime}$.
The universal property of the Hilbert scheme gives a map
$U^{\prime} \ra H_{4,-3a-b}$ whose image is precisely the family
$F_{1}$ of unions $W \cup_{3P} L$.
In particular, $F_{1}$ is irreducible.

The structure of the map $U^{\prime} \rightarrow H$ shows that
$U^{\prime}$ has dimension $\binom{a+b+5}{3}-a-2b+2$.
On the other hand, if $W$ is a triple line arising in the construction
above, then $\ideal{W}$ is $(a+b+2)$-regular (see \cite[2.4]{nthree}),
hence $\dim H^{0}(\ideal{W}(a+b+2))= \binom{a+b+5}{3}-6a-4b-9$.
Subtracting this redundancy shows
that the family has dimension $5a+2b+11=11-2g-a$.

The proof of statement 2 goes through via the same outline. The main
differences are as follows. The family $H \subset H_{3,-a-1}$ is
now the family of unions $Z \cup_{P} L$, which is irreducible of
dimension $10+2a$ by \cite[3.2(b)]{nthree}. In the definition \ref{v}
of $V$, $I_{Y}^{3}$ is replaced by $I_{Y}^{2}$ and in the map
\ref{mp} the exponents of $y$ should be reduced by one. To see that
$U^{\prime}$ is nonempty, we can use the same triple line $W$ as in
the proof above, but instead use the line $L$ given by $\{x=z=0\}$.
The remaining modifications are clear.
\end{proof}

\br\label{bigfibres}{\em Since the family of triple lines $W$ is
irreducible of dimension $10+5a+2b$ \cite[2.6]{nthree}, we expect
that the natural map $\{W \cup_{3P} L\} \stackrel{F}{\rightarrow}
\{W\}$ which forgets the line $L$ has generically one dimensional
fibres. However, there are triple lines $W$ for which the fibre
$F^{-1}(W)$ has larger dimension. For example, the triple lines
constructed in characteristic $p>0$ by Hartshorne
\cite[2.3]{genus} have a two-dimensional family of triple tangent
lines. \em}\er

\bp\label{WLclosure} With the notation above, we have
\begin{enumerate}
       \item $F_{2} \subset {\overline{F_{1}}}$ if
             $\ds 0 \leq a < \frac{-g-1}{3}$ and
             is otherwise empty.
       \item $F_{4} \subset {\overline{F_{3}}}$ if
             $\ds 0 \leq a < \frac{-g-1}{3}$ and
             is otherwise empty.
\end{enumerate}
\ep

\begin{proof}
Let $W_{0} \cup_{2P} L$ be a curve in the family $F_{2}$, so that the
underlying double line $Z \subset W_{0}$ satsifies $\Length (Z \cap L) = 2$.
If $W_{0}$ has type $(a,-3a-g-1)$ and support $Y$; we may write
$I_{Y}=(x,y)$, $I_{Z} =((x,y)^{2}, xg-yf)$, $I_{L} = (x,z)$ and
$I_{W_{0}} = (I_{Y} I_{Z}, h_{0})$ in suitable coordinates \cite[2.3]{nthree}.
If $K \subset H^{0}(\ideal{Z}(-2a-g+2))$ is the vector subspace
considered in the proof above, then $zh_{0} \in K$.
Fixing a member $(Z \cup L, h) \in U^{\prime}$ as above, the deformation
$(1-t)zh_{0} + th$ gives a map
$\Aone \stackrel{\psi}{\rightarrow} K$, which yields
$\psi^{-1}(U^{\prime}) \rightarrow H_{4,g}$.
This extends to a map
${\overline \psi}:T = \psi^{-1}(U^{\prime}) \cup \{0\} \ra {\rm Hilb_{4}^{g}}$
into the full Hilbert scheme: by construction, it's clear that the
limit curve ${\overline \psi}(0)$ contains $W_{0} \cup L$. Since this
curve has genus $g$, it is equal to ${\overline \psi}(0)$, completing
the proof. The limit of the triple lines $W_{t}$ is the triple line
$W_{0}$ along with an embedded point, which is conveniently covered up
by the line $L$. Statement 2 is similar.
\end{proof}

\bc\label{extend} The closure of family $F_{1,a}$ in $H_{4,g}$ contains
extremal curves.
\ec

\begin{proof}
Following \cite[3.6]{nthree}, the family of ideals $I_{t}$ below
give a deformation from a triple line $W$ of type $(a,b)$ to an
extremal curve of the same arithmetic genus. $$I_{t}=((x,y)^{3},
(x,y)(xz^{a+1}-tyw^{a+1}),z^{b} t^{2} (xz^{a+1}-tyw^{a+1}) - x^{2}
w^{a+b})$$ We simply observe that the line $L=\{x=w=0\}$ is triple
tangent to the triple line $W_{t}$ defined by $I_{t}$ for all $t
\neq 0$ and that this same line is a triple tangent to the limit
extremal curve having ideal
$$I_{0}=(x^{2},xy,y^{3},xz^{3a+b+3}-y^{2}w^{3a+b+2}).$$
\end{proof}

\br\label{closure}{\em The closure of the family $F_{3}$ above forms an
irreducible component of the Hilbert scheme (Theorem \ref{comps}) with
one exception. When $a=0$ and $g \leq -2$, a curve
$C_{0} = W_{0} \cup_{2P} L \in F_{1}$
is a flat limit of curves $C_{t} = Z_{t} \cup_{2P} L \cup L_{t}$,
where $Z_{t}$ is a double line and $L_{t}$ is a line disjoint from $L$.
This is not surprising in view of \cite[Proposition 3.3]{nthree},
which says that $W$ is the limit of unions $Z_{t} \cup L_{t}$.
% The deformation is similar to those already given in this section.

To see this, let $Z$ be the underlying double line and
$Y = \supp W$. First we use \cite[Proposition 2.6]{nthree} to write
$I_{W}=((x,y)^{3}, xq, yq, h=pq-ax^{2}-bxy-cy^{2})$ where
$q=xz-yw$ may be taken to be the equation of a {\it smooth}
quadric surface $Q$ by \cite[Remark 1.5]{nthree} and $I_{L}=(x,l)$
for some linear form $l$. Note that $L$ is not tangent to $Q$ at
$P$ because $L$ meets the underlying double line $Z$ in a reduced
point.

If $L_{t}=\{x+wt = y+zt = 0\}$ for $t \in \Aone$, the family
$D_{t}=L_{t} \cup Y \cup_{P} L$ gives a flat family of extremal
curves whose limit is $D_{0}=Z \cup L$. In considering the total
family $D \subset \Pthree \times \Aone \stackrel{\pi}{\rightarrow}
\Aone$, we see by Grauert's theorem that
$\pi_{*}(\ideal{D}(-g+1))$ is locally free on $\Aone$ (the
extremal curves have constant cohomology) and hence globally free.
Since $\ideal{D_{1}}(3)$ is generated by global sections and
$D_{1}$ is reduced of embedding dimension $\leq 2$, $D_{1}$ lies
on smooth surfaces of degrees $\geq 3$. In particular, we can find
a section $s_{t}$ (yielding a corresponding surface $S_{t}$) such
that $S_{1}$ is a smooth surface containing $D_{1}$ and $s_{0} = l
h$.

Now consider the family $C_{t}=S_{t} \cap (L \cup Y^{\{2\}} \cup
L_{t})$. For general $t \neq 0$, $C_{t}$ is the disjoint union of
$L_{t}$ and $Z_{t} \cup_{P} L$, where $Z_{t}$ is a double line of
genus $g$. Let $D_{0}$ be the flat limit in the (full) Hilbert
scheme. It is clear that $D_{0} = L \cup W_{0}$ $ +$ possible
embedded points, with $W_{0}$ a triple line. Since the limit of
$L_{t} \cup Y$ is $Z$, it is clear that $Z \subset W_{0}$ and
hence $W_{0}$ has type $(0,b^{\prime})$ for some $b^{\prime} \geq
0$. Since $l h \in I_{W_{0}}$ and the plane $\{l=0\}$ meets
$W_{0}$ properly, we have that $h \in I_{W_{0}}$ and it follows
(from \cite[Proposition 2.6]{nthree}) that $W_{0} \subset W$.
Since these are Cohen-Macaulay triple lines, we conclude that
$W_{0} = W$ and hence $D_{0} = L \cup W$ (there are no embedded
points because $p_{a}(C_{t})=p_{a}(L \cup W)=g$).

\em}\er

\section{The Hilbert Schemes $H_{4,g}$} \label{summary}

In this section we prove the main results of the paper. The first
of these is Theorem \ref{comps} that describes the irreducible
components of the Hilbert schemes $H_{4,g}$. The second is the fact
\ref{hilbconnected} that $H_{4,g}$ is connected. The cases when
$g \geq 0$ are well known and described in the introduction. We begin
with the case $g=-1$, since it has a somewhat different statement due to the
existence of more reduced curves.

\bp \label{p4-1}
The Hilbert scheme $H_{4,-1}$ is connected and has $3$
irreducible components:
\begin{enumerate}
         \item
         The $17$-dimensional family of extremal curves.
         \item
         The $16$-dimensional closure of the family of subextremal
         curves. The general member is the disjoint union of two conics.
         \item
         The $16$-dimensional closure of the family whose general member is
         a disjoint union of a twisted cubic and a line.
\end{enumerate}
\ep
\begin{proof}
       If a curve $C \in H_{4,-1}$ is not extremal, then its spectrum is
       necessarily $\{0^{2},1^{2}\}$, in which case $C \in H_{2}$ or
       $H_{3}$ by Proposition \ref{d3g0-c}. The families $H_{2}$ and $H_{3}$
       have general members as described in Remark \ref{descriptions} and
       meet because both contain thick 4-lines. Finally, $H_{2}$ meets the
       family of extremal curves by \cite{nolletp} or \cite{2h}.
\end{proof}

In the following theorem the letter $L$ denotes a line and the
symbol $\dot{\cup}$ the disjoint union of two curves. $L \cup_{nP} C$
denotes the schematic union of a line $L$ and a curve $C$,
whose intersection is the divisor $nP$ on $L$.
\bt\label{comps}
The irreducible components of the Hilbert schemes $H_{4,g}$ for $g
\leq -2$ are those listed in the following table \vspace{.2cm}
%%%%%%%%%tab
\begin{center} \begin{em}
\begin{tabular}{|c|c|c|c|}
\hline Label   & General Curve & Dimension & Restrictions \\

\hline $G_{1}$ &
\begin{tabular}{c}
$D \cup Z$ \\ $D$ smooth conic, $\deg (Z) = 2$ \\  $g(Z) = g-3$,
$\text{length} (D \cap Z) = 4$
\end{tabular}
& $15-2g$  &none \\
    \hline

& $L_{1} \cup_{2P} Z \cup_{2Q} L_{2}$
     &  &\\
$G_{2}$ &$L_{1} \cap L_2 = \emptyset$  & $13-2g$& none\\ & $\deg
(Z) = 2$,  $g(Z) = g-2$ & & \\ \hline

& $D \cup_{2P} Z$
    &  &  \\
$G_{3}$ & $D$ smooth conic&$13-2g$&none \\ & $\deg (Z) = 2$, $g(Z)
= g-1$  & & \\ \hline

&&& \\ $G_{4}$ & general thick $4$-line &  $9-3g$ & $g \leq -3$ \\
&&& \\ \hline

&&& \\ $G_{5}$ & double conic  & $13-2g$ & none \\ & && \\ \hline

$G_{6}$ &
\begin{tabular}{c}
$Z \cup_{2P} L_{1} \dot{\cup} L_{2}$ \\ $\deg (Z) = 2$,  $g(Z) =
g$
\end{tabular}
&
\begin{tabular}{c} \\
$11-2g$  \\  \\
\end{tabular}
&$g \leq -3$
\\ \hline
$G_{7,a}$ &
\begin{tabular}{c}
    $W \cup_{3P} L$
\\
$W$ quasiprimitive $3$-line
\\
type of $W = (a,-3a-g)$
\end{tabular}
&$11-2g-a$ &
\begin{tabular}{c}
$g \leq -3$ \\ $0 < a \leq  \frac{-g}{3}$
\end{tabular}
\\
\hline

$G_{8,a}$ &
\begin{tabular}{c}
    $W \cup_{2P} L$
\\
$W$ quasiprimitive $3$-line
\\
type of $W = (a,-1-3a-g)$
%\\ $\text{length} (W_2 \cap L) =1$
\end{tabular}
&$10-2g-a$ &
\begin{tabular}{c}
$g \leq -6$ \\ $0 < a <  \frac{-g-1}{3}$
\end{tabular}
\\
\hline

$G_{9,a}$ &
\begin{tabular}{c}
$W \dot{\cup} L$
    \\
$W$ quasiprimitive $3$-line
\\
    type of $W = (a,-3-3a-g)$
\end{tabular}
& $8-2g-a$ &
\begin{tabular}{c}
$g \leq -6$ \\ $0 < a \leq  \frac{-g-3}{3}$
\end{tabular} \\
\hline

&  $D_1 \dot{\cup} D_2$ &  & \\ $G_{10,m}$ & $\deg( D_{1}) =2$,
$g (D_1) = -m$ & d(-m) + d(g+m+1) &$0 \leq m \leq \frac{-g-1}{2}$
\\ &$\deg (D_2) = 2$,  $g(D_{2})=g+m+1$& & \\ \hline

$G_{11,a,b}$&
\begin{tabular}{c}
Quasiprimitive $4$-line \\
    type $(a,b,-6a-b-g-3)$
\end{tabular}
& \begin{tabular}{c} $7-2g-3a =$ \\   $9a+2b+2c+13$
\end{tabular}
& \begin{tabular}{c}$g \leq -9$\\ $0 < a \leq \frac{-g-3}{6}$\\ $0
\leq b \leq \frac{-6a-g-3}{2}$
\end{tabular}
\\
\hline

\end{tabular}  \end{em}
\end{center}
\vspace{.2in} where $d(-m) = \dim \text{H}_{2,-m}$ equals $5+2m$
if $m > 1$, and $8$ if $m=0$ or $m=1$.
\et
\begin{proof}
In the table $G_{j}$ denotes the closure in the Hilbert scheme
$H_{4,g}$ of the set of curves described in the corresponding row.
The outline of the proof is as follows. First we show  the families
listed in the table are irreducible of the stated dimension. Then we
show there is no inclusion relation among them. Finally, we prove every curve
of degree $4$ and genus $g \leq -2$ belongs to one of these families.
We will see the restrictions given in
the table are necessary to ensure that a given family exists and
is not contained in another family of the list.

The family $G_1$ consists of extremal curves, and is an
irreducible component of the Hilbert scheme of dimension $15-2g$
by \cite[4.3]{extremal}. $G_2$, $G_3$ and $G_4$ are the closures
of the families $H_2$, $H_3$ and $H_{n(g)}$
of Proposition~\ref{d3g0-c}, and are therefore irreducible of the stated
dimension. Note that $G_2$ contains the subextremal curves, and $G_4$
consists of thick $4$-lines.
The closure $G_5$ of the family of double conics of genus $g$ is
irreducible of dimension $13-2g$ because the family of such curves
in a {\it fixed} double plane $2H$ is irreducible of dimension
$10-2g$ by \cite[2.1 and 4.3]{2h}. $G_{7,a}$ and $G_{8,a}$ are the
closures of families $F_{1,a}$ and $F_{3,a}$ from Proposition~\ref{WL}, and
are therefore irreducible of the stated dimension.
The irreducibility of $G_{11,a,b}$ is proven in Proposition~\ref{qp}.

Finally, families $G_{6}$, $G_{9,a}$ and $G_{10,m}$ are irreducible
components of $H_{4,g}$ because the curves defining them are
disjoint unions of curves that are general in their respective
Hilbert schemes. The dimensions of these families can be computed
out of \cite[1.6,3.4,3.5]{nthree}.

This shows all families in the statement are irreducible, and
we now prove there are no inclusion relation among them.
We have just seen that $G_1$, $G_{6}$, $G_{9,a}$
and $G_{10,m}$ are irreducible components of $H_{4,g}$,  so
certainly none of them can be contained in any other family of the
list.

The families $G_{2}$, $G_{3}$ and $G_{5}$ could only be contained
in $G_{1}$ or $G_{4}$  by reason of dimension. However, none of
them is contained in $G_{1}$ by semicontinuity, and none of them
is contained in $G_4$ because curves in $G_4$ are supported on a
single line.

In the table, the only families of dimension larger than that of
$G_4$ consist of curves $C$ with $h^0 \idealc (2) \neq 0$. These
cannot specialize to a general thick $4$-line $T$ because $h^0
\ideal{T} (2) = 0$ for $g \leq -3$ by Proposition~\ref{d3g0-c}.
Thus $G_{4}$ is not contained in any of the other families. Note
however that $G_{3}$ contains all thick $4$-lines when $g=-2$.

We now treat the case of $G_{7,a}$ and $G_{8,a}$. The general
curves in $G_{7,a}$ and $G_{8,a}$ are not supported on a line, so
they can't be contained  in $G_4$ or $G_{11,a,b}$. $G_{7,a}$ and
$G_{8,a}$ are not contained in  $G_{1}$, $G_{2}$ $G_3$ or $G_{5}$
by semicontinuity - curves in the latter families have larger
speciality by Proposition~\ref{d3g0-c}.

$G_{7,a}$ and $G_{8,a}$ are not contained in $G_{6}$ or $G_{10,m}$
because when two lines collapse the resulting double line has
genus $\geq -1$, hence $a \leq 0$ contradicting the restriction
imposed.

If $G_{7,a^{\prime}} \subset  G_{9,a}$, then in considering the
underlying double line as in Remark \ref{qplinespec} we see that
$a^{\prime} \leq a$, which in turn implies that $\dim
G_{7,a^{\prime}}= 11-2g-a^{\prime} > 8-2g-a = \dim G_{9,a}$, a
contradiction. Similarly $G_{8,a^{\prime}} \not\subset G_{9,a}$.

It remains to show that neither $G_{7,a}$ and $G_{8,a}$ contains
the other. There can be no containment $G_{7,a^{\prime}} \subset
G_{8,a}$, because then $a^{\prime} \leq a$ by Remark
\ref{qplinespec} and hence $\dim G_{8,a}= 10-2g-a <
11-2g-a^{\prime} = \dim G_{7,a}$, a contradiction. Now suppose
that $G_{8,a^{\prime}} \subset G_{7,a}$. Remark \ref{qplinespec}
tells us again that $a^{\prime} \leq a$, and since $\dim
G_{8,a^{\prime}} < \dim G_{7,a}$ we conclude that $a^{\prime} =
a$. In particular, the limit of the underlying family of double
lines $Z$ has no embedded points. This is not possible because the
limit double line meets $L$ in one point while the general member
meets $L$ in a double point.

Finally, $G_{11,a,b}$ cannot be contained in any of the families
$G_{j}$ with $j \leq 10$ by semicontinuity: indeed, since $a > 0$,
every $4$-line $C$ in $G_{11,a,b}$ satisfies $h^{1} \coo_{C} (-2)
= 1$, while for any other curve $D \in H_{4,g}$ we have $h^{1}
\coo_{D} (-2) \geq 2$.
%because for $a > 0$, the curves contain no degree
%two subcurve lying on four independent surfaces of degree two
%(this holds for the rest of the families).
On the other hand, there are no containments among the families
$G_{11,a,b}$: if $G_{11,a,b} \subset
G_{11,a^{\prime},b^{\prime}}$, then by Remark \ref{qplinespec} we
would have $a \leq a^{\prime}$, while $$7 -2g -3a^{\prime} = \dim
G_{11,a^{\prime},b^{\prime}} > \dim G_{11,a,b} = 7 -2g -3a$$ would
yield $a^{\prime} < a$, a contradiction.

To finish the proof, we still have to show our families cover the Hilbert
scheme. Let $C \in H_{4,g}$ have support $B = C_{red}$.

{\vskip .10in} \noindent{ \bf Case 1:} $\deg B = 4$ {\vskip .10in}

\noindent Here $C=B$ is reduced, and  all reduced curves of degree
$4$ satisfy $g \geq -1$ with the following two exceptions: either
(a) $C$ is the disjoint union of a conic (possibly degenerate) and
two lines, when $g =-2$ and $C \in G_{10,0}$ or (b) $C$ is the
disjoint union of four lines, $g=-3$ and $C \in G_{10,1}$.

\vskip .10in \noindent{ \bf Case 2:} $\deg B = 3$ \vskip .10in

\noindent In this case $C = Z \cup D$, where $Z$ is a double line
with support $L$ and $D$ is a reduced curve of degree $2$. In
particular, $B=L \cup D$.

First suppose that $D$ is planar and let $l = \Length (D \cap L)$.
If $l=0$, then $C \in G_{10,0}$. If $l=1$, then $g(B)=0$ and $C$
belongs to one of the families $G_{2}$, $G_{3}$ or $G_{4}$ by
Proposition \ref{d3g0-c}. If $l=2$, then $B$ is planar and hence
$C$ is extremal by \cite[2.2]{subextremal}.

The other possibility is that $D=L_{1} \cup L_{2}$ is a disjoint
union of lines. If $D$ does not meet $Z$, then $C \in G_{10,1}$.
If $Z$ meets $L_1$ but not $L_2$, then $Z \cup L_{1}$ is a
specialization of a double line meeting $L_{1}$ in a double point
by \cite[3.2]{nthree}, so $C$ lies in $G_{6}$. If $Z$ meets both
$L_1$ and $L_2$, then $g(B)=0$ and $C$ again belongs to one of the
families $G_{2}$, $G_{3}$ or $G_{4}$ by Proposition \ref{d3g0-c}.

\vskip .10in

\noindent{ \bf Case 3:} $\deg B = 2$ \vskip .10in

\noindent If $B$ is a smooth conic, then $C$  belongs to
$G_{5}$, so we may assume that $B=L \cup L^{\prime}$ is a union of
two lines. If $C$ is a union of two double lines, then either (a)
the lines are disjoint and $C \in G_{10,m}$ for $m=-\max
\{g(Z_{1}), g(Z_{2})\}$ or (b) the lines meet and hence $C$ is
contained in a double plane; in this case $C \in G_{5}$ by
\cite[8.1 and 8.2]{2h}.

The remaining possibility is that $C=W \cup L$, where $L$ is a
line and $W$ is a triple line which is quasi-primitive because $g
(W) \leq -1$ (the only thick triple line has genus $0$). If $W$
has type $(a,b)$,
let $Z$ be the underlying double line % with $g(Z)=-a-1$
and set $l=\Length (W \cap L)$.

If $l=3$, then necessarily $\Length (Z \cap L)=2$. If $a = -1$,
then $Z \cup L$ is planar and $C$ is extremal. If $a=0$, then $Z
\cup L \in H_{3,0}$ and $C$ belongs to one of the families
$G_{2}$, $G_{3}$ or $G_{4}$ by Proposition \ref{d3g0-c}. If $a >
0$, then $C$ belongs to $G_{7,a}$.

Suppose $l=2$. If $\Length (Z \cap L)=2$ and $a=-1$ or $0$, we
argue as in the case $l=3$. If $\Length (Z \cap L)=2$ and $a > 0$,
then $C \in G_{7,a}$ by Proposition~\ref{WLclosure}. Thus we may assume
$Z \cap L=P$ a {\it reduced} point. If $a=-1$, then $Z \cup L \in
H_{3,0}$ and $C \in G_{2} \cup G_{3}$ by Proposition \ref{d3g0-c}.
If $a = 0$, then $C \in G_{6}$ by Remark \ref{closure}. If $a >
0$, then $C \in G_{8,a}$.

Suppose $l=1$. If $a > 0$,  then $C$ is in $G_{8,a}$ by
Proposition~\ref{WLclosure}. If $a = 0$, then $C$ belongs to $G_{6}$ by
Proposition~\ref{WLclosure} and Remark \ref{closure}.

If $a=-1$, $Z \cup L \in H_{3,0}$ and $C
\in G_{2} \cup G_{3}$ by Proposition~\ref{d3g0-c}.

If $l=0$, then $C = W \cup L$ with $W$ quasi-primitive of type
$(a,-3-3a-g)$. If $a=-1$, then $W$ is extremal and $C \in G_{6}$
by \cite[3.2]{nthree}. If $a=0$, then $C \in G_{5,0}$ by
\cite[3.3]{nthree}, while if $a > 0$, then $C \in G_{9,a}$.

\vskip .10in \noindent{ \bf Case 4:} $\deg B = 1$ \vskip .10in

\noindent If $C$ is thick, then $C \in G_{4}$. If $C$ is a
quasiprimitive $4$-line, then $C$ has type $(a,b,c)$ for some
integers $a \geq -1$ and $c \geq b \geq 0$. If $a=-1$, then the
underlying double line $Z$ is planar and $C$ lies in a double
plane, hence $C \in G_{5}$ by \cite[8.1 and 8.2]{2h}. If $a=0$,
then $C$ is in $G_{10,m}$ for $m=-1-b$ by Proposition
\ref{disjointdoublelines}. Finally, if $a > 0$, then $C \in
G_{11,a,b}$.

\vskip .10in

\end{proof}

\bex\label{collapse}{\em As the restrictions in the statement
imply, some of these components do not show up if the genus $g$ is
not small enough.
%For example, if $g=-2$,
%the curves in $G_{10,a,b}$ simply don't exist,
%since $a > 0$ forces the genus of the $4$-line to be $\leq -6$.
%The family $G_{10,0}$ is the closure of disjoint unions
%of a conic and two lines. Since $H_{3,-1}$ is irreducible
%\cite[3.1]{nthree}, both $G_{5}$ and $G_{9,a}$ are
%contained in $G_{10,0}$. The family $G_{4}$ is in the closure
%of $G_{2}$ and $G_{3}$ because $j=2$ or $3$ in Proposition~\ref{d3g0-c}.
%
%Now we consider the families $G_{7,a}$ and $G_{8,a}$. Since $a > 0
%\Rightarrow p_{a}(W) \leq -5 \Rightarrow p_{a}(W \cup_{3P} L) \leq -3$,
%the family $G_{7,a}$ is empty here. Similar genus considerations
%reduce us to the looking at $G_{8,0}$ so that $W$ has type $(0,0)$
%and genus $-1$. Applying \cite[3.3]{nthree}, we find that
%$G_{8,0} \subset {\overline G_{5}}$. Thus only five
%components remain.
%
%Summarizing,
%$H_{4,-2}$ is connected and has $5$ irreducible components:
For example, $H_{4,-2}$  has only $5$ irreducible components:
\begin{enumerate}
\item
The $19$-dimensional family of extremal curves $G_1$, whose
general member is the union of a conic and a double line of genus
$-5$ meeting in a scheme of length $4$.
\item
The $17$-dimensional family $G_2$, whose general member, the union
$L_{1} \cup_{2P} Z \cup_{2Q} L_{2}$ of lines $L_{i}$ and a double
line $Z$ of genus $-4$, is a subextremal curve.
\item
The $17$-dimensional family $G_3$ whose general member is the
union of a double line of genus $-3$ and a conic, meeting in a
double point.
\item
The $17$-dimensional family $G_5$ whose general member is a double
conic.
\item
The $16$-dimensional family  $G_{10,0}$ of curves whose general
member is the disjoint union of a conic and two lines.
\end{enumerate}
\em}\eex

We can now prove $H_{4,g}$ is connected:
\bt\label{hilbconnected}
The Hilbert scheme $H_{4,g}$
is connected whenever nonempty.
\et

\begin{proof}
We may assume $g \leq -2$, and it suffices to show that all the
irreducible components can be connected to the component $G_1$ of
extremal curves. The families $G_{2}$ and $G_{3}$ meet the family
$G_{4}$ of thick $4$-lines by Proposition \ref{d3g0-c}, and
$G_{2}$  meets $G_{1}$ by \cite{phsamir}, \cite{2h}, or
\cite{nolletp}. In particular, thick $4$-lines belong to the
connected component of extremal curves, and it is enough to show
that all other components can be connected to either $G_1$ or
$G_4$. $G_{5}$ meets $G_{1}$ by \cite[5.1 and 8.2]{2h}. That
$G_{6}$ meets $G_{1}$ follows immediately  from \cite[2.1]{hconn}
and~\cite{nthree}. $G_{7,a}$ meets $G_{1}$ by
Corollary~\ref{extend} and  $G_{8,a}$ meets $G_{1}$ by
Proposition~\ref{WLclosure} and \cite[2.5]{hconn}. One can connect
$G_{9,a}$ to $G_{1}$ by applying \cite[3.8]{nthree} and
\cite[2.1]{hconn}. The families $G_{11,a,b}$ meet $G_{4}$ by
Proposition \ref{thintothick}. Applying Proposition
\ref{disjointdoublelines}, we see $G_{10,m}$
         contains $G_{11,0,m-1}$ for $m > 0$, and hence meets $G_{4}$ as
         well.
         Lastly, we consider $G_{10,0}$.
         By definition
         this family contains curves $C = Z \cup L_{1} \cup_{P} L_{2}$
         where $Z$ has degree two and $L_{i}$ are meeting lines. Since
         $Z$ is extremal, $C$ specializes to an extremal
         curve  $E \in G_1$ by \cite[2.1 and 2.5]{hconn}.

\end{proof}

As an application of our results, we can now give a counterexample
to a conjecture of
A{\"\i}t-Amrane and Perrin \cite{aitperrin}. The conjecture regards
the following question, which has been a recurring theme of
this paper:

\bq{\em
Let $X$ and $X_{0}$ be two irreducible families of curves in
$H_{d,g}$ {\it having constant cohomology}. Under what conditions do we
have a nonempty intersection ${\overline X} \cap X_{0} \neq \emptyset$
in $H_{d,g}$?
\em}\eq

We have been lucky in that whenever we suspected the existence
of such a deformation, we could actually prove it. In general, this
question is difficult. A first necessary condition is
provided by semicontinuity \cite[III, $\S 12$]{AG}: If
$\tau^{i}(n) = h^{i} (\pso, \ideal{C}(n))$
for $C \in X$ and
$\tau^{i}_{0}(n) = h^{i} (\pso, \ideal{C_{0}}(n))$ for
$C_{0} \in X_{0}$, then whenever ${\overline X} \cap X_{0} \neq \emptyset$
we must have $\tau^{i}(n) \leq \tau^{i}_{0}(n)$ for all $i$ and $n$
(we write $\tau \leq \tau_{0}$ for short). This condition is not
sufficient, even when $X$ (resp. $X_{0}$) is an irreducible component
of the Hilbert scheme $H_{\tau}$ (resp. $H_{\tau_{0}}$) of curves with
fixed cohomology. This has been shown by a recent example of
A{\"\i}t-Amrane and Perrin \cite{aitperrin}.

A more subtle necessary condition is afforded by the Rao modules of
the curves. Let $A$ be a discrete valuation ring with fraction field
$K$ and residue field $k$ and let
$\mathcal{C} \subset \pso_{A}$ be a family of locally
Cohen-Macaulay curves over $A$. Then
the Rao module $M_{C_{K}}$ of the generic curve
is a flat deformation of a subquotient $M$ of the Rao module
$M_{C_{k}}$ of the special curve (\cite[Proposition 5.9]{hmdp3},
\cite[$\S 4.2.2$]{phsamir}, \cite[Proposition 13]{aitperrin}). This
means that there are submodules $M_{1} \subset J \subset M_{C_{k}}$ such
that $M = J/M_{1}$ and $M$ is a flat deformation of $M_{C_{K}}$. Moreover,
\begin{equation} \label{ait}
\dim (M_1)_n = h^0 \ideal{C_{k}} (n) -
h^0 \ideal{C_{K}} (n) \,\,\, \mbox{and} \,\,\,
\dim (M_{C_{k}}/J)_n = h^2 \ideal{C_{k}} (n) -
h^2 \ideal{C_{K}} (n)
\end{equation}

In view of this result, the following conjecture of
A{\"\i}t-Amrane and Perrin is natural:

\bconj[\cite{aitperrin}, Conjecture 14]
Let $X$ and $X_0$ be irreducible components of $H_{\tau}$ and
$H_{\tau_0}$ respectively.
Suppose that
\begin{enumerate}
\item
$\tau \leq \tau_{0}$
%$\tau^{i} (n) \leq \tau_{0}^{i} (n)$ for all $ i \geq 0$ and $n \in \Z$;
\item
The Rao module of the generic curve $C_{\xi}$ of $X$ is a flat deformation of a
subquotient of the Rao module of a curve $C_0$ in $X_{0}$ and that
the numerical
conditions (\ref{ait}) hold.
\end{enumerate}
Then ${\overline X} \cap X_0 \neq \emptyset$ in $H_{d,g}$.
\econj

As it turns out, the conditions of the conjecture are still
insufficient, as we note in the following example.

\bex\label{perrin}{\em
For $g \leq -3$, let $X \subset H_{4,g}$ denote the irreducible family
$H_{n(g)}$ from Proposition \ref{d3g0-c}. The closure ${\overline X}$
is precisely the family of thick 4-lines. Let $X_{0}=G_1$ denote
irreducible family consisting of extremal curves. We claim that
${\overline X} \cap X_{0} = \emptyset$. Indeed, an extremal curve cannot
specialize to a thick 4-line because this would violate semicontinuity, while
a thick 4-line has everywhere embedding dimension $3$, and so cannot specialize
to an extremal curve that has generic embedding dimension $2$~\cite{extremal}.
On the other hand, we will now show that the conditions of the conjecture hold.

First we compare the Rao modules. Let $C$ be a general thick
4-line with support $L$ and set $S = S_{L}$. Proposition \ref{d3g0-c} shows
that $M_{C} \cong S/(a,b,c)(-g-1)$ where $a,b,c$ are general forms of
degree $-g+1$ in $S$. Choose a linear form $l \in S$ so that
$(a,l)$ is a regular sequence and a form $f \in (b,c)$ of degree
$-g+2$ so that $(a,f)$ is a regular sequence.
We consider the extremal Koszul module $M = S/(a,lf)$. Since the
multiplication $S/(a) \stackrel{\cdot l}{\rrr} S/(a)$ is injective and the
image of the submodule $({\overline f})$ is $({\overline {lf}})$, we
see that the submodule $J = {\overline l}M$ is isomorphic to $S/(a,f)$.
Since $f \in (b,c)$ by choice, $S/(a,b,c)$ is a quotient of $J$ by
$M_{1}=(lb,lc)$. If $E$ is an extremal curve corresponding to $M$,
then $\deg E = 4$, $p_{a}(E)=g$ and we have just shown that $M_{C}$ is a
subquotient of $M_{E}$. It is clear that $\dim (M/J)_{n}=1$ for
$g \leq n \leq 0$ and zero otherwise; this is seen to be precisely
$h^2 \ideal{E} (n) - h^2 \ideal{C} (n)$ in comparing the spectra of
these curves. Finally, since the Euler characteristics of $\ideal{E}$
and $\ideal{C}$ are the same, the exact sequences relating the Rao
modules shows that
$\dim (M_{1})_{n} = h^0 \ideal{E} (n) - h^0 \ideal{C} (n)$. In
particular, the semicontinuity conditions are immediate.

\em}\eex

\bibliographystyle{alpha}

\end{document}